\RequirePackage{etoolbox}
\csdef{input@path}{{style/}{graphics/}}
\documentclass[aos,MSNbibl,seceqn,nameyear,dvips]{arximspdf}
\usepackage{graphicx}
\doi{10.1214/14-AOS1279}
\volume{43}
\issue{1}
\pubyear{2015}
\firstpage{352}
\lastpage{381}
\docsubty{FLA}

\makeatletter
\newcommand{\lleft}{\left}
\newcommand{\rrvert}{\vert}
\newcommand{\rright}{\right}
\newcommand{\llvert}{\vert}
\renewcommand{\mid}{|}
\newtheorem{theorem}{Theorem}[section]
\newtheorem{prop}[theorem]{Proposition}
\newproclaim{defn}[theorem]{Definition}
\newproclaim{rem}[theorem]{Remark}
\newproclaim{ex}{Example}
\newcommand{\iint}{\int\!\!\int}
\newcommand{\eqref}[1]{(\ref{#1})}
\def\bbeta{\bolds{\beta}}
\def\bX{\mathbf{X}}
\def\risk{\operatorname{Risk}}
\def\tglrt{\mathsf{T}_\mathrm{GLRT}}
\def\tglrtnew{\mathsf{T}_\mathrm{GLRT}}
\def\thc{\mathsf{T}_\mathrm{HC}}
\def\thcnew{\mathsf{T}_\mathrm{HC}}
\def\I{\mathcal{I}}
\def\bnu{\bolds{\nu}}
\def\BC{\mathrm{BC}}
\makeatother

\begin{document}
\begin{frontmatter}

\title{Hypothesis testing for high-dimensional sparse binary regression}
\runtitle{Hypothesis testing}

\begin{aug}
\author[A]{\fnms{Rajarshi}~\snm{Mukherjee}\corref{}\ead[label=e1]{rmukherj@stanford.edu}},
\author[B]{\fnms{Natesh S.}~\snm{Pillai}\ead[label=e2]{pillai@fas.harvard.edu}}
\and
\author[C]{\fnms{Xihong}~\snm{Lin}\ead[label=e3]{xlin@hsph.harvard.edu}}
\runauthor{R. Mukherjee, N.~S. Pillai and X. Lin}
\affiliation{Stanford University, Harvard University and Harvard University}
\address[A]{R. Mukherjee\\
Department of Statistics\\
Stanford University\\
Sequoia Hall\\
390 Serra Mall\\
Stanford, California 94305-4065\\
USA\\
\printead{e1}}
\address[B]{N.~S. Pillai\\
Department of Statistics\\
Harvard University\\
1 Oxford Street\\
Cambridge, Massachusetts 01880\\
USA\\
\printead{e2}}
\address[C]{X. Lin\\
Department of Biostatistics\\
Harvard University\\
655 Huntington Avenue\\
SPH2, 4th Floor\\
Boston, Massachusetts 02115\\
USA\\
\printead{e3}}
\end{aug}

%
\received{\smonth{6} \syear{2014}}
%
\revised{\smonth{10} \syear{2014}}

%
\begin{abstract}
In this paper, we study the detection boundary for minimax hypothesis
testing in the context of high-dimensional, sparse binary regression
models. Motivated by genetic sequencing association studies for rare
variant effects, we investigate the complexity of the hypothesis
testing problem when the design matrix is sparse. We observe a new
phenomenon in the behavior of detection boundary which does not occur
in the case of Gaussian linear regression. We derive the detection
boundary as a function of two components: a design matrix sparsity
index and signal strength, each of which is a function of the sparsity
of the alternative. For any alternative, if the design matrix sparsity
index is too high, any test is asymptotically powerless irrespective of
the magnitude of signal strength. For binary design matrices with the
sparsity index that is not too high, our results are parallel to those
in the Gaussian case. In this context, we derive detection boundaries
for both dense and sparse regimes. For the dense regime, we show that
the generalized likelihood ratio is rate optimal; for the sparse
regime, we propose an extended Higher Criticism Test and show it is
rate optimal and sharp. We illustrate the finite sample properties of
the theoretical results using simulation studies.
\end{abstract}

%
\end{frontmatter}

\section{Introduction}\label{sec1}
The problem of testing for the association between a set of covariates
and a response is of fundamental statistical interest. In the context
of testing for a linear relationship of covariates with a continuous
response, R.~A.~Fisher introduced analysis of variance (ANOVA) in the
1920s, which is still widely used in the present day. In recent years,
finding the detection boundary of various testing problems has gained
substantial popularity. A fruitful way of finding the detection
boundary is to study the minimax error of testing and obtain a
threshold of signal strength under which all testing procedures in the
concerned problem are useless. For Gaussian linear models, this has
been extensively studied by \citet{Candes} and \citet{Ingster5}; these
works were inspired by the previous work on hypothesis testing in
various contexts, such as sparse normal mixtures [\citet
{Jin1,Cai1}], Gaussian sequence models [\citet{Ingster4}] and
correlated multivariate normal problems [\citet{Jin2}]. However,
very little work has been done on detection boundaries in generalized
linear models for discrete outcomes.

In this paper, we study the detection boundary for hypothesis testing
in the context of high-dimensional, sparse binary regression models.
Motivated by case--control sequencing association studies for detecting
the effects of rare variants on disease risk [\citet
{tang2013large,lee2014}],
we are interested in the complexity of the hypothesis testing problem
when the design matrix is sparse. { Specifically, sequencing studies
allow sequencing massive genetic variants in candidate genes or across
the whole genome. A rapidly increasing number of sequencing association
studies have been conducted, such as the 1000 Genome Project
[\citet{10002012integrated}]
and the NHLBI Exome Sequencing Project [\citet{fu2013analysis}].
It is of substantial interest to study rare variant effects on diseases
case--control candidate gene and whole genome sequencing association
studies. A~major challenge in analysis of sequencing data is that a
vast majority of variants across the genome are rare variants [\citet{10002012integrated}
(Figure~2b), \citet
{fu2013analysis} (Figure 1a), \citet{nelson2012abundance} (Figure 1c)].
For a review of analysis of data of sequencing association studies, see
\citet{lee2014}. }

\begin{figure}

\includegraphics{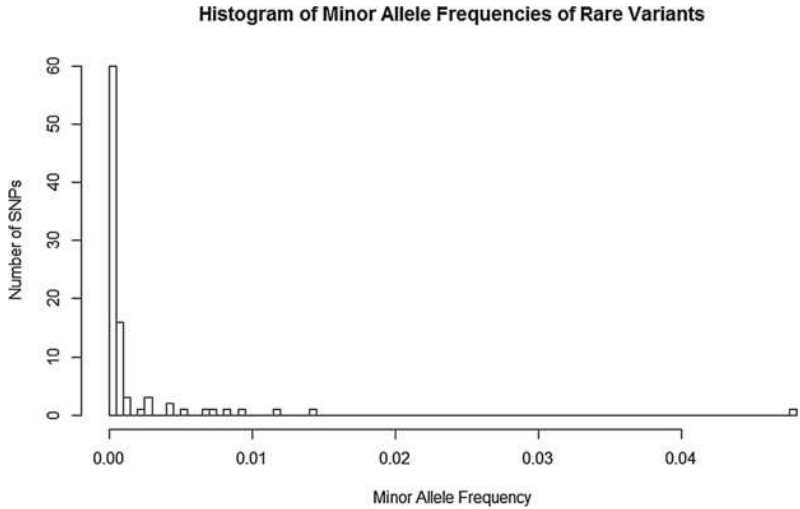}

\caption{The histogram of minor allele frequencies of uncommon/rare
variants (MAF${}\leq{}$5\%) in the Dallas Heart study data.} \label{DHSMAF}
\end{figure}

For example, in the Dallas Heart candidate gene sequencing study
[\citet{dallas}], 3476 individuals were sequenced in the region
consisting of three genes ANGPTL3, ANGPTL4 and ANGPTL5. The goal of
study was to test the effects of these genes on the risk of
hypertriglyceridemia. A total of 93 genetic variants were observed in
these genes. Each variant took values 0, 1, 2, which represents the
number of minor alleles in a genetic variant. About half of the
variants were singletons, {that is}, they were observed in only
one person; 92 variants have the minor allele frequencies${}<5\%$.
The design matrix is hence very sparse, with a vast majority of its
columns having $<$5\% nonzero values (1 or 2), and the proportion of total nonzero
elements in the design matrix being $<$2.5\%. It is expected only a
small number of variants might be associated with hypertriglyceridemia.
The presence of the sparse design matrix and sparse signals for binary
outcomes results in substantial challenges in testing the association
of these genes and hypertriglyceridemia. Figure~\ref{DHSMAF} provides
the histogram of rare variants with minor allele frequencies less than
$5\%$. 

%

Suppose there are $n$ samples of binary outcomes, $p$ covariates for
each. Consider a binary regression model linking the outcomes to the
covariates. We are interested in testing a global null hypothesis that
the regression coefficients are all zero and the alternative is sparse
with $k$ signals, where $k=p^{1-\alpha}$ and $\alpha\in[0,1]$. For
binary regression models, we observe a new phenomenon in the behavior
of detection boundaries which does not occur in the Gaussian framework,
as explained below.

The main contribution of our paper is to derive the detection boundary
for binary regression models as a function of two components: a design
matrix sparsity index and signal strength, each of which is a function
of the sparsity of the alternative, {that is}, $\alpha$.
Throughout the paper, we will call the first component as ``design
matrix sparsity index.'' This is unlike the results in Gaussian linear
regression which has a one component detection boundary, namely the
necessary signal strength. In the Gaussian linear model framework,
\citet{Candes} and \citet{Ingster5} show that if the design matrix
satisfies certain ``low coherence conditions,'' then it is possible to
detect the presence of a signal in a global sense, provided the signal
strength exceeds a certain threshold. In contrast, our results suggest
that for binary regression problems, the difficulty of the problem is
also determined by the design matrix sparsity index. In this paper, we
explore two key implications of this phenomenon which are outlined
below.

First, if the design matrix sparsity index is too high, we show that no
signal can be detected irrespective of its strength. In Section~\ref
{secSDM}, we provide sufficient conditions on the design matrix
sparsity index which yield such nondetectability problems. Such
conditions on the design matrix sparsity index corresponds to the first
component of the detection boundary. \citeauthor{plan1} (\citeyear{plan1,plan2}) discussed a
difficulty in inference similar to that of ours, for design matrices
with binary entries in the context of 1-bit compressive sensing and
sparse logistic models. Our results in Section~\ref{secSDM} pertain
to sparse design matrices with \emph{arbitrary entries}, which are not
necessarily orthogonal. We give a few examples of design matrices which
satisfy our criteria for nondetectability. These include block diagonal
matrices and banded matrices.  

Second, for design matrices with binary entries and with low
correlation among the columns, we are able to characterize both
components of the detection boundary. In particular, if the design
matrix sparsity index, the first component of the detection boundary,
is above a specified threshold, no signal is detectable irrespective of
strength.
Once the design matrix sparsity index is below the same threshold, we
also obtain the optimal thresholds with respect to the second component
of the detection boundary, {that is}, the minimum signal
strength required for successful detection. In this regime, our results
parallel the theory of detection boundary in Gaussian linear
regression. We also provide relevant tests to attain the optimal
detection boundaries. In the sparse regime $(\alpha>\frac{1}{2})$,
our results are sharp and rate adaptive in terms of the signal strength
component of the detection boundary. Moreover, we observe a phase
transition in both components of the detection boundary depending on
the sparsity $(\alpha)$ of the alternative. To the best of our
knowledge, this is the first work optimally characterizing a two
component detection boundary in global testing problems against sparse
alternatives in binary regression. 

To illustrate further, we contrast our results with the existing
literature. In the case of a balanced one-way ANOVA type design matrix
with each treatment having $r$ independent replicates, for Gaussian
linear models, \citet{Candes} show that the detection boundary is given
by $O(\frac{p^{1/4}}{\sqrt{kr}})$ in the dense regime
($\alpha<\frac{1}{2}$) and
equals $\sqrt{\frac{2\rho_{\mathrm{linear}}^{*}(\alpha)\log (p)}{r}}$ in the sparse regime
$\alpha>\frac{1}{2}$, where
%
%
\begin{equation}
\label{eqndetboundLM} \rho_\mathrm{linear}^{*}(\alpha) = \cases{
\alpha-\frac{1}{2}, &\quad if $ \frac
{1}{2} < \alpha<
\frac{3}{4}$,
\cr
(1-\sqrt{1-\alpha})^2, & \quad if $\alpha\geq
\frac{3}{4}$}
\end{equation}
and $\rho_\mathrm{linear}^{*}(\alpha)$ matches the detection
boundary in \citet{Jin1} in the normal mixture problem. For given
sparsity of the alternative, the detection boundary depends a single
function of $r$.

For binary regression, we show that the detection boundary is
drastically different and 
depends on two functions of $r$: a design matrix sparsity index and
signal strength under the alternative hypothesis for a given regime. In
particular, define the \textit{design matrix sparsity index} of a design
matrix as $1/r$. For $r=1$, every test is powerless irrespective of the
signal sparsity and the signal strength under the alternative hypothesis.
When $r>1$, the behavior of the detection boundary can be categorized
into three situations. In\vspace*{1pt} the \emph{dense} regime where $r>1$ and
$\alpha\leq\frac{1}{2}$, the detection boundary matches that of the
Gaussian case up to rates and the usual Generalized Likelihood Ratio
Test achieves the detection boundary. In the \emph{sparse} regime,
{that is}, when $\alpha> \frac{1}{2}$, the detection boundary
behaves differently for $r \ll\log(p)$ and $r \gg\log(p)$. For
$\alpha> \frac{1}{2}$ and $r \ll\log(p)$, a new phenomenon that
does not exist in the Gaussian case arises: all tests are
asymptotically powerless irrespective of how strong the signal strength
is in the alternative. For\vspace*{1pt} $\alpha> \frac{1}{2}$ and $r \gg\log
(p)$, our results are identical to the Gaussian case, up to a constant
factor accounting for the Fisher information.
In this regime, we construct a version of the Higher Criticism Test and
show that this test achieves the lower bound.
We use the strong embedding theorem [\citet{kmt}] to obtain sharp
detection boundary. Noting that this problem can also be cast as a test
of homogeneity among~$p$ binomial populations with contamination in $k$
of them.
Hence, roughly speaking, the two component detection boundary in this
binary problem setting equals $[1,O(\frac{p^{1/4}}{\sqrt
{kr}}))$ in the dense regime and $(O(\frac{1}{\log(p)}),O(\sqrt
{\frac{\log(p)}{r}})]$ in the sparse regime, where the first
component represents the design matrix sparsity index, which is of the
order of $1/r$,
and the second component indicates the order of signal strength.
Successful detection requires both components to be above the
component-specific detection boundaries.

Borrowing ideas from orthogonal designs, we further obtain analogous
results for general binary design matrices which are sparse and have
weak correlation among columns, mimicking design matrices often
observed in sequencing association studies. For such general binary
designs, we are able to completely characterize the two component
detection boundary in both dense and sparse regimes. Our versions of
Generalized Likelihood Ratio Test and the Higher Criticism Test
continue to attain the optimal detection boundaries in dense and sparse
regimes, respectively. Similar to orthogonal designs, our results are
sharp in the sparse regime and we once again obtain optimal phase
transition in the two component detection boundary depending on the
sparsity $(\alpha)$ of the alternative. Our results show that under
certain low correlation structures, the problem essentially behaves as
an orthogonal problem.

The rest of the paper is organized as follows. We first formally
introduce the model in Section~\ref{secprelim} and discuss general
strategies. Here, we also provide a set of notation to be used
throughout the paper. In Section~\ref{secSDM}, we study the
nondetectability for sparse design matrices with arbitrary entries. In
Section~\ref{secDD}, we formally introduce a class of designs for
which we derive the sharp detection boundaries, namely, one-way ANOVA
designs and weakly correlated binary designs. Section~\ref{sectests}
introduces the Generalized Likelihood Ratio Test (GLRT) and the Higher
Criticism Test in our designs, which will be used in subsequent
sections to attain the sharp detection boundaries in two different
regimes of sparsity. In Section~\ref{secSADE}, we first analyze the
one-way ANOVA designs and derive the sharp detection boundary in
different sparsity regimes. In Section~\ref{secWA}, we derive the
sharp detection boundary in different sparsity regimes for weakly
correlated binary designs. Section~\ref{secNE} presents simulation
studies which validate our theoretical results. Finally, we collect all
the technical proofs in the supplementary material [\citet
{MukherjeePillaiLin2014}].

\section{Preliminaries}\label{sec2}\label{secprelim}
Suppose there are $n$ binary observations $y_i \in\{0,1\}$, for $1
\leq i \leq n$, with covariates $\mathbf{x}_i = (x_{i1},\ldots
,x_{ip})^t$. The design matrix with rows $\mathbf{x}_i^t$ is denoted
by $\mathbf{X}$.
Set ${\mathbf y} = (y_1, y_2,\ldots, y_n)^t$.
The conditional distribution of $y_i$ given $\mathbf{x}_i$ is given by
%
%
\begin{equation}
\label{eqnmodmain} \mathbb{P}(y_i=1\mid\mathbf{x}_i,
\bbeta)=\theta \bigl(\mathbf{x}_i^t\bbeta \bigr),
\end{equation}
where $\bbeta=(\beta_1,\ldots,\beta_p)^t \in\mathbb{R}^p$ is an
unknown $p$-dimensional vector of regression coefficients.
Henceforth, we will assume that $\theta$ is an arbitrary distribution
function that is symmetric around $0$, {that is},
%
%
\begin{equation}
\label{eqntheta} \theta(z)+\theta(-z)=1\qquad\mbox{for all } z \in \mathbb{R}.
\end{equation}
For some of the results, we will also require certain smoothness
assumptions on $\theta(\cdot)$ which we will state when and where
required. Examples of such $\theta(\cdot)$ include logistic and
normal distributions which, respectively, correspond to logistic and
probit regression models. 

Let $M(\bbeta)=\sum_{j=1}^p {I(\beta_j \neq0)}$ and let $R_k^p=\{
\bbeta\in\mathbb{R}^p\dvtx M(\bbeta) = k\}$. For some \mbox{$A > 0$}, we are
interested in testing the global null hypothesis
%
%
\begin{equation}
\label{eqnhyp} H_0\dvtx \bbeta=0\quad\mbox{vs}\quad H_1
\dvtx \bbeta\in \Theta^A_k= \biggl\{ \bbeta\in\bigcup
_{k' \geq k}R_{k'}^p\dvtx \min \bigl\{
\llvert\beta _j\rrvert\dvtx \beta_j \neq0 \bigr\} \geq A
\biggr\}.\hspace*{-20pt}
\end{equation}
Set $k=p^{1-\alpha}$ with $\alpha\in(0,1]$. We note that these types
of alternatives have been considered by \citet{Candes}, referred to as
the ``\textit{Sparse Fixed Effects Model}'' or SFEM. In particular,
under the alternative, $\bbeta$ has at least $k$ nonzero coefficients
exceeding $A$ in absolute values. Alternatives corresponding to $\alpha
\leq\frac{1}{2}$ belong to the \textit{dense regime} and those
corresponding to $\alpha> \frac{1}{2}$ belong to the \textit{sparse
regime}. We will denote by $\pi$ a prior distribution on $\Theta
^A_k\subset\mathbb{R}^p$. Throughout we will refer to $A$ as the
signal strength corresponding to the alternative in equation~\eqref{eqnhyp}.

We first recall a few familiar concepts from statistical decision
theory. Let a test be a measurable function of the data taking values
in $\{0,1\}$. The Bayes risk of a test $\mathsf{T} = \mathsf
{T}(\mathbf{X}, {\mathbf y})$ for testing $H_0\dvtx \bbeta= 0$ versus
$H_1\dvtx \bbeta\sim\pi$ when $H_0$ and
$H_1$ occur with the same probability, is defined as the sum of its
probability of type I error (false
positives) and its average probability of type II error (missed detection):
\[
\risk_{\pi}(\mathsf{T}):=\mathbb{P}_0(\mathsf{T}=1)+\pi
\bigl[\mathbb {P}_{\bbeta}(\mathsf{T}=0)\bigr],
\]
where $\mathbb{P}_{\bbeta}$ denotes the probability distribution of
${\mathbf y}$ under model \eqref{eqnmodmain} and $\pi[\cdot]$ is the
expectation with respect to the prior $\pi$. 
We study the asymptotic properties of the binary regression model
\eqref{eqnmodmain} in the high-dimensional regime, {that is},
with $p \rightarrow\infty$ and $n=n(p) \rightarrow\infty$ and a
sequence of priors $\{\pi_p\}$. Adopting the terminology from \citet
{Candes}, we say that a sequence of tests $\{\mathsf{T}_{n,p}\}$ is
\textit{asymptotically powerful}\break if $\lim_{p \rightarrow\infty
}\operatorname{Risk}_{\pi_p}(\mathsf{T}_{n,p})=0$, and\vspace*{1pt} it is \textit
{asymptotically powerless} if\break $\liminf_{p \rightarrow\infty}\mathrm
{Risk}_{\pi_p}(\mathsf{T}_{n,p})\geq1$. When no prior is specified,
the risk is understood to be the worst case risk or the minimax risk
defined as
\[
\risk(\mathsf{T}):=\mathbb{P}_0(\mathsf{T}=1)+\max
_{\bbeta\in
\Theta^A_k} \bigl[\mathbb{P}_{\bbeta}(\mathsf{T}=0) \bigr].
\]
%

The detection boundary of the testing problem \eqref{eqnhyp} is the
demarcation of signal strength $A$ which determines whether all tests
are asymptotically powerless (we call this lower bound of the problem)
or there exists some test which is asymptotically powerful (we call
this the upper bound of the problem).

To understand the minimax risk, set
\[
d(\mathcal{P}_0,\mathcal{P}_1)=\inf\bigl\{\llvert P-Q
\rrvert _1\dvtx P \in\mathcal{P}_0, Q \in
\mathcal{P}_1\bigr\},
\]
where $\mathcal{P}_0, \mathcal{P}_1$ are two families of probability
measures and $\llvert P-Q\rrvert _1=\break  \sup_{B}\llvert
P(A)-Q(A)\rrvert $, with $B$ being a
Borel set in $\mathbb{R}^n$, denotes the
total-variation norm. Then for any test $\mathsf{T}$, we have
[\citet{wald1}]
\[
\risk(\mathsf{T}) \geq1-\tfrac{1}{2}\,d\bigl(\mathbb{P}_0,
\operatorname {conv}_{\bbeta\in\Theta^A_k}(\mathbb{P}_{\bbeta})\bigr),
\]
where $\operatorname{conv}$ denotes the convex hull. However, $d(\mathbb
{P}_0,\operatorname{conv}_{\bbeta\in\Theta^A_k}(\mathbb{P}_{\bbeta}))$
is difficult to calculate. But it is easy to see that for any test
$\mathsf{T}$ and any prior $\pi$, one has $\operatorname{Risk}(\mathsf
{T}) \geq\operatorname{Risk}_{\pi}(\mathsf{T})$. So in order to prove
that a sequence of tests is asymptotically powerful, it suffices to
bound from above the worst-case risk $\operatorname{Risk}(\mathsf{T})$.
Similarly, in order to show that all tests are asymptotically
powerless, it suffices to work with an appropriate prior to make
calculations easier and bound the corresponding risk from below for any
test $\mathsf{T}$.

It is worth noting that, for any prior $\pi$ on the set of $k$-sparse
vectors in $\mathbb{R}^p$ and for any test $\mathsf{T}$, we have
\[
\risk_{\pi}(\mathsf{T}) \geq1-\tfrac{1}{2}\mathbb{E}_{0}
\llvert L_{\pi
}-1\rrvert \geq1-\tfrac{1}{2}\sqrt{
\mathbb{E}_0\bigl(L_{\pi}^2\bigr)-1},
\]
where $L_{\pi}$ is the $\pi$-integrated likelihood ratio and $\mathbb
{E}_0$ denotes the expectation under~$H_0$. For the model \eqref{eqnmodmain}, we have
%
%
\begin{equation}
\label{eqnintlike} L_{\pi}=2^n\int{\prod
_{i=1}^n{ \biggl(\frac{\theta(\mathbf{x}_i^t
\bbeta)}{\theta(-\mathbf{x}_i^t \bbeta)}
\biggr)^{y_i}\theta \bigl(-\mathbf{x}_i^t \bbeta
\bigr)}\,d\pi(\bbeta)}.
\end{equation}
Hence, in order to assess the lower bound for the risk, it suffices to
bound from above $\mathbb{E}_0(L_{\pi}^2)$. By Fubini's theorem, for
fixed design matrix $\mathbf{X}$, we have
%
%
\begin{eqnarray}\label{eqnintlikesq}
\qquad && \mathbb{E}_0 \bigl(L_{\pi}^2
\bigr)
\nonumber\\[-8pt]\\[-12pt]
&&\qquad = 2^n \iint\prod_{j=1}^n
{ \bigl[\theta \bigl(\mathbf{x}_i^t \bbeta \bigr)\theta
\bigl(\mathbf{x}_i^t \bbeta' \bigr)+\theta
\bigl(-\mathbf{x}_i^t \bbeta \bigr)\theta \bigl(-
\mathbf{x}_i^t \bbeta' \bigr) \bigr]}\,d
\pi( \bbeta)\,d\pi \bigl(\bbeta' \bigr),\nonumber
\end{eqnarray}
where $\bbeta,\bbeta' \sim\pi$ are independent.
In the rest of the paper, all of our analysis is based on studying
$\mathbb{E}_0(L_{\pi}^2)$ carefully for the prior distribution $\pi$
chosen below.

In the context of finding an appropriate test matching the lower bound,
by the Neyman--Pearson lemma, the test which rejects when $L_{\pi}>1$
is the most powerful Bayes test and has risk equal to $1-\frac
{1}{2}\mathbb{E}_{0}\llvert L_{\pi}-1\rrvert $. However, this test requires
knowledge of the sparsity index $\alpha$ and is also computationally
intensive. Hence, we will construct tests which do not require
knowledge of $\alpha$ and are computationally much less cumbersome.

Ideally, one seeks least favorable priors, {that is}, those
priors for which the minimum Bayes risk equals the minimax risk.
Inspired by \citet{baraud1}, we choose $\pi$ to be uniform over all
$k$ sparse subsets of $\mathbb{R}^p$ with signal strength either $A$
or $-A$. 

\subsection{Notation}\label{sec2.1}
We provide a brief summary of notation used in the paper. For two
sequences of real numbers $a_p$ and $b_p$, we say $a_p \ll b_p$ or
$a_p=o(b_p)$, when $\limsup_{p \rightarrow\infty} \frac
{a_p}{b_p} \rightarrow0$ and we\vspace*{1pt} say $a_p \lesssim b_p$ or $a_p=O(b_p)$
if\break $\limsup_{p \rightarrow\infty} \frac{a_p}{b_p}<\infty$.
The indicator function of a set $B$ will be denoted by $\I(B)$.

We take $\pi$ to be uniform over all $k$ sparse subsets of $\mathbb
{R}^p$ with signal strength either $A$ or $-A$. Let $M(k, p)$ be the
collection of all subsets of $\{1,\ldots, p\}$ of size $k$. For each
$m \in M(k, p)$, let $\xi^m =(\xi_j)_{j \in m}$ be a sequence of
independent Rademacher random variables taking values in $\{+1,-1\}$
with equal probability. Given $A > 0$ for testing \eqref{eqnhyp}, a
realization from the prior distribution $\pi$ on $\mathbb{R}^p$ can
be expressed as
\[
\beta_{\xi,m} =\sum_{j \in m}{A\xi_j
e_j},
\]
where $(e_j)_{j=1}^p$ is the\vspace*{1pt} canonical basis of $\mathbb{R}^p$ and $m$
is uniformly chosen from $M(k,p)$. Since, the alternative in \eqref
{eqnhyp} allows both positive and negative directions of signal
strength $\beta_j$, we call it a two-sided alternative. On the
contrary, when we are given the extra information in \eqref{eqnhyp}
that the $\beta_j$'s have the same sign, then we call the alternative
a one-sided alternative. A realization from a prior distribution over
one-sided $k$ sparse alternatives can be expressed as $\sum_{j \in
m}{A\xi e_j}$, where $\xi$ is a single Rademacher random variable.

For any distribution $\pi'$ on $M(k,p)$, by $\textbf{support}(\pi
')$ we mean the smallest set $I':=\{M\dvtx M \in M(k,p)\}$ such that
$\pi'(I')=1$. For any distribution $\pi^*$ over $M(k,p)$, we
say that another distribution $\pi_0$ over $M(k,p)$ is equivalent to
$\pi^*$ (denoted by $\pi_0 \sim\pi^*$) if $\pi_0$ is uniform on
its support and
\[
\pi^* \bigl(M \notin\textbf{support}(\pi_0) \bigr)=o(1).
\]
By the support of a vector $v \in\mathbb{R}^p$, we mean the set $\{j
\in\{1,\ldots,p\}\dvtx v_j \neq0\}$; the vector $v$ is $Q$-sparse if the
support of $v$ has at most $Q$ elements. For $i=1,\ldots,n$, we will
denote the support of the $i$th row of $\bX$ by $S_i:=\{
j\dvtx \bX_{i,j}\neq0\}\subset\{1,\ldots,p\}$. Let $\BC^l$ denote the
set of all functions whose $l$th derivative is continuous
and bounded over $\mathbb{R}$. By $\theta(\cdot) \in\BC^l(0)$, we
mean that the $l$th derivative of $\theta(\cdot)$ is
continuous and bounded in a neighborhood of $0$. Finally, by saying
that a sequence measurable map $\chi_{n,p}(\mathbf{}y,\bX)$ of the
data is tight, we mean that it is stochastically bounded as $n,p
\rightarrow\infty$. 

\section{Sparse design matrices and nondetectability of signals}\label{sec3}
\label{secSDM}
In this section, we study the effects of sparsity structures of the
design matrix $\bX$ on the detection of signals. Our key results in
Theorem \ref{theoremnodetectcondition} below provide a sufficient
condition on the sparsity structure of the $\bX$ which renders all
tests asymptotically powerless in the sparse regime irrespective of
signal strength $A$. This result for nondetectability is quite general
and { are satisfied by different classes of sparse design matrices
as we discuss below.}
We verify the hypothesis of Theorem \ref{theoremnodetectcondition}
in a few instances where certain global detection problems can be
extremely difficult.\vadjust{\goodbreak}

Let $\pi_0 \sim\pi$ and $R_{\pi_0}$ denote the support of $\pi_0$.
For a sequence of positive integers~$\sigma_p$, we say that $j_1,j_2
\in\{1,\ldots,p\}$ are $``\sigma_p$-mutually\vspace*{1pt} close'' if $\llvert
j_1-j_2\rrvert
\leq\sigma_p$. For an $m_1$ from $\pi_0$ and $N\geq0$, let
${R}^N_{m_1}(\sigma_p)$ denote the set of all $\{l_1,\ldots,l_k\} \in
R_{\pi_0}$ such that there are exactly $N$ elements $``\sigma
_p$-mutually close'' with members of $m_1$.

\begin{theorem}\label{theoremnodetectcondition}
Let $k=p^{1-\alpha}$ with $\alpha> \frac{1}{2}$. Let $\pi_0 \sim
\pi$ and $\{\sigma_p\}$ be a sequence of positive integers with
$\sigma_p \ll p^{\varepsilon}$ for all $\varepsilon>0$. Let $m_1$ be drawn
from $\pi_0$. Suppose that for all $N=0,\ldots,k$ and every $m_2$
drawn from $\pi_0$ with $m_2 \in{R}^N_{m_1}(\sigma_p)$, the
following holds for some sequence $\delta_p > 0$:
%
%
\begin{equation}
\label{eqnsparsecondition} \sum_{i=1}^n \bigl\{\I
\bigl(\min \bigl\{\llvert m_1 \cap S_i\rrvert,\llvert
m_2 \cap S_i\rrvert \bigr\}>0 \bigr) \bigr\} \leq N
\delta_{p},
\end{equation}
where $S_i$ is defined in the last paragraph of Section~\ref{sec2}.
Then if
$\delta_{p} \ll\log(p)$, all tests are asymptotically powerless.
\end{theorem}

An intuitive explanation of Theorem \ref{theoremnodetectcondition}
is as follows. If the support of $\bbeta$ under the alternative does
not intersect the support of a row of the design matrix $\bX$, the
observation corresponding to that particular row does not provide any
information about the alternative hypothesis. If randomly selected
draws from $M(k,p)$ fail to intersect with the support of most of the
rows, as quantified by equation \eqref{eqnsparsecondition}, then all
tests will be asymptotically powerless irrespective of the signal
strength in the alternative. 
In the Gaussian linear regression, the effect of a similar situation is
different. We provide an intuitive explanation for a special case in
Section~\ref{sec6}. Also intuitively, the quantity $\frac{1}{\delta_p}$ in
Theorem \ref{theoremnodetectcondition} is a candidate for the design
matrix sparsity index of $\bX$. This is because if $\frac{1}{\delta
_p}$ is too large, as quantified by $\frac{1}{\delta_p} \gg\frac
{1}{\log(p)}$, then all tests are asymptotically powerless in the
sparse regime irrespective of the signal strength.
Now we provide a few examples where condition \eqref
{eqnsparsecondition} can be verified to hold for appropriate parameters.

\begin{ex}[(Block structure)] 
Suppose that, up to permutation of rows, $\bX$ can be partitioned into
a block diagonal matrix consisting of $\mathbf{G}^{(1)},\ldots,\break 
\mathbf{G}^{(M)}$ and a matrix $\mathbf{G}$ as follows:
%
%
\begin{equation}
\label{eqnblockdiagonalmatrix} \mathbf{X}= \lleft[ %
\begin{array}
{ccccc} \mathbf{G}^{(1)}_{c_1 \times d_1} & & & &
\\
& \ddots& & &
\\
& & \mathbf{G}^{(j)}_{c_j \times d_j} & &
\\
& & & \ddots&
\\
& & & & \mathbf{G}^{(M)}_{c_M \times d_M}
\\
\hline& & \raisebox{-4pt} {{\mbox{${\mathbf{G}}_{\tilde{c}
\times
p}$}}} & &
\\[-5pt]
&&&& \end{array} %
\rright] \in\mathbb{R}^{n \times p},
\end{equation}
where $\tilde{c} = n - \sum_{j=1}^M c_j$.
The matrices $\mathbf{G},\mathbf{G}^{(1)},\ldots, \mathbf{G}^{(M)}$
are \emph{arbitrary matrices} of specified dimensions. Let $c^*=\max_{1\leq j\leq M}c_j$ and $l^*=\max_{1\leq j\leq
M}d_j$. Indeed $c^*, l^*$ and the structure of $\mathbf{G}$ decide the
sparsity of the design matrix $\bX$. In Theorem \ref
{theoremblockdesign} below, we provide necessary conditions on $c^*,
l^*$ and $\mathbf{G}$ which dictate the validity of condition \eqref
{eqnsparsecondition}, and hence renders all tests asymptotically
powerless irrespective of signal strength.

%
\begin{figure}

\includegraphics{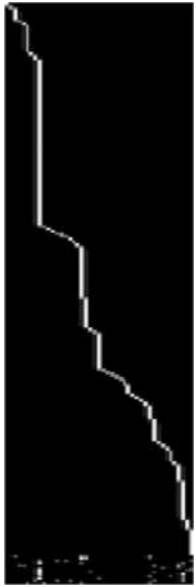}

\caption{Heat map of the genotype matrix $\bX$ of the Dallas Heart
study data after a suitable rearrangement of subject indices, after
removing the single common variant. The nonzero entries of the genotype
matrix that represent mutations are colored white, while the zero
entries that represent no mutation are colored in black.}
\label{DHS}
\end{figure}

 Design matrices in sequencing association studies for rare variants
generally have this structure. Figure~\ref{DHS} shows a heat map of
the genotype matrix of the subjects in the Dallas Heart study after a
suitable rearrangement of subject indices, after removing the single
common variant. It shows that the genotype matrix has the same
structure as $\bX$ described above. Specially, it can be partitioned
into two parts. The top part of the matrix is an orthogonal block
diagonal structure and the bottom part is a nonorthogonal sparse matrix
which corresponds to $\mathbf{G}$.
\end{ex}

%

\begin{theorem}\label{theoremblockdesign}
Assume\vspace*{1pt} that the matrix $\bX$ is of the form given by \eqref
{eqnblockdiagonalmatrix}. Let $k=p^{1-\alpha}$ with $\alpha> \frac
{1}{2}$ and\vspace*{1pt} suppose that $\llvert \bigcup_{i > {n^*}} S_i\rrvert \ll p$
where $n^*=\sum_{j=1}^M c_j$. Let 
${l^*} \ll{p^{\varepsilon}}$ for all $\varepsilon>0$. If $c^* \ll\log
{p}$, then condition \eqref{eqnsparsecondition} holds, and thus all
tests are asymptotically powerless.
\end{theorem}


In Theorem \ref{theoremblockdesign}, the condition $\llvert \bigcup_{i>n^*}S_i\rrvert \ll p$ is an assumption on the structure of $\mathbf{G}$
which restricts the locations of nonzero elements of $\mathbf{G}$.
This condition on~$\mathbf{G}$ is not tight and can be much relaxed
provided one assumes further structures on~$\mathbf{G}$. In fact, this
implies that asymptotically the bulk of the information about the
alternatives comes from the block diagonal part of $\bX$ and the
information from $\mathbf{G}$ is asymptotically negligible.


{Further, intuitively, $\frac{1}{c^{*}}$ is the candidate for the
design matrix sparsity index. Since if $\frac{1}{c^{*}}$ is too high,
as quantified by $\frac{1}{c^{*}} \gg\frac{1}{\log(p)}$, then all
tests are asymptotically powerless in the sparse regime. It is natural
to ask about the situation when the design matrix sparsity index is
below the specified\vspace*{2pt} threshold of $\frac{1}{\log(p)}$, {that
is}, $c^* \gg\log(p)$.} To this end, it is possible to analyze the
necessary and sufficient conditions on the signal strength $A$
dictating asymptotic detectability in problem \eqref{eqnhyp} when
$c^* \gg\log(p)$ for $\bX$ in \eqref{eqnblockdiagonalmatrix} but
possibly with $\llvert \bigcup_{i >n^*}S_i\rrvert \gg p$.
In Section~\ref{secWA}, we provide an answer to this question when
$\bX$ has binary entries.



\begin{ex}[(Banded matrix)]
Suppose $\bX$ has the following banded structure, possibly after a
permutation of its rows.
Suppose there exists $l_2>l_1$ such that for $i=1,\ldots,n$,
$X_{i,j}=0$ for $j<i-l_1$ or $j>i+l_2$. { Further, let $\llvert \bigcup_{i > n} S_i\rrvert \ll p$.} 
Note that this allows design matrices $\bX$ which can be partitioned
into a banded matrix of band-width $l_2-l_1$ and an arbitrary design
matrix with sparsity restrictions as specified by $\llvert \bigcup_{i
> n} S_i\rrvert \ll p$.
\end{ex}

%
%
\begin{theorem}
\label{theorembandedmatrix}
Let\vspace*{1pt} $k=p^{1-\alpha}$ with $\alpha> \frac{1}{2}$. Suppose $\bX$ is a
banded design matrix as described above.
Suppose that $l_2-l_1 \ll\log(p)$. Then condition \eqref
{eqnsparsecondition} holds and thus all tests are asymptotically powerless.
\end{theorem}
\section{Design matrices}\label{sec4} \label{secDD}In Section~\ref
{secSDM}, we
provided conditions on $\bX$ under which all tests are asymptotically
powerless irrespective of signal strength $A$. To complement those
results, the subsequent sections will be devoted toward analyzing
situations when $\bX$ is not pathologically sparse, and hence one can
expect to study \mbox{nontrivial} conditions on the signal strength $A$ that
determine the complexity in~\eqref{eqnhyp}. In this section, we
introduce certain design matrices with binary entries motivated by
sequencing association studies. In subsequent sections, we will derive
the detection boundary for binary regression models with these design
matrices.

In order to introduce the design matrices we wish to study, we need
some notation. 
Set $\Omega^*=\{i\dvtx \llvert S_i\rrvert =1\}$. For $j=1,\ldots,p$,
let $\Omega^*_j=\{
i\in\Omega^*\dvtx S_i=\{j\}\}$ with $r_j=\llvert \{i\in\Omega^*\dvtx
S_i=\{j\}\}
\rrvert $. Let $r^*=\max_{1 \leq j \leq p}r_j $ and $r_*=\min_{1 \leq
j \leq p}r_j$. Also, let $n^*=\sum_{j=1}^p r_j$ and
$n_{*}=n-n^*$. In words, for each $j$, $\Omega^*_j$ is the collection
of individuals with only one nonzero informative covariate appearing as
the $j$th covariate and $r_j$ is the number of such individuals.

A binary design matrix, as described above, is orthogonal if and only
if all of its rows have at most one nonzero element. Hence, up to a
permutation of rows, any binary design matrix can be potentially
partitioned as a one-way ANOVA type design and an arbitrary matrix. In
particular, up to a permutation of rows, any binary design matrix is
equivalent to equation \eqref{eqnblockdiagonalmatrix} where each
$\mathbf{G}^{(j)}_{r_j \times1}=(1,\ldots,1)^t$, $c_j=r_j$,
$d_j=1$, $c^*=r^*$, $l^*=1$, $\tilde{c}=n_*$ and $\mathbf{G}$ is an
arbitrary matrix with binary entries. Keeping this in mind, we have the
following definitions.
%

\begin{defn}\label{def1}
\label{wa}
A design matrix $\bX$ is defined as a Weakly Correlated Design with
parameters $(n^*,n_*,r^*,r_*,Q_{n,p},\gamma_{n,p})$ if the following
conditions hold:

\begin{longlist}[(C3)]
\item[(C1)] The design matrix $\bX_{n \times p}$ has binary entries;

\item[(C2)] $\llvert S_i\rrvert \leq Q_{n,p}$ for all $i=1,\ldots
,n$, for some
sequence $Q_{n,p}$;

\item[(C3)] $\frac{n_{*} Q_{n,p}^2}{r^*} \ll\gamma_{n,p}$ for some
sequence $\gamma_{n,p} \rightarrow\infty$.
\end{longlist}
\end{defn}

As a special case of the above definition, we have the following definition.
%

\begin{defn}\label{def2}
\label{sa}
A design matrix $\bX$ is called an ANOVA design with parameter $r$,
and denoted by $\bX\in\mathrm{ANOVA}(r)$, if it is a Weakly Correlated Design with $r_*=r^*=r$ and $n_*=0$.
\end{defn}

A few comments are in order for the above set of assumptions in
Definitions \ref{def1} and \ref{def2}. The motivation for condition
\textup{(C1)} comes from genetic association studies assuming a dominant
model. As our proofs will suggest, this can be easily relaxed, allowing
the elements of $\bX$ to be uniformly bounded above and below.
Condition~\textup{(C2)} imposes sparsity on $\bX$. Finally,
since the part of $\bX$ without $\mathbf{G}$ is exactly orthogonal,
condition \textup{(C3)} restricts the deviation of $\bX$ from exact
orthogonality. In particular, if the size of $\mathbf{G}$ is
``not too large'' compared to the orthogonal part of $\bX$, as we will
quantify later, then the behavior of the detection problem is similar
to the one with an exactly orthogonal design. In essence, this captures
low correlation designs suitable for binary regression with ideas
similar to low coherence designs as imposed by \citet{Candes} for
Gaussian linear regression. 

Because of the presence of $\mathbf{G}$, Weakly Correlated
Designs in Definition \ref{def1} allow for correlated binary design
matrices with sparse structures. However, condition \textup{(C3)}
restricts the size of $\mathbf{G}$ (numerator) compared to the
orthogonal part (denominator) by a factor of $\gamma_{n,p}$.
Intuitively, this implies low correlation structures in~$\bX$.
The condition \textup{(C3)} restricts the effect of $\mathbf{G}$ on the
correlation structures of $\bX$ by not allowing too many rows compared
to the size of the orthogonal part of $\bX$.
It is easy to see that when $n_*Q_p \ll p$, then since $\llvert \bigcup_{i \notin\Omega^*}S_i\rrvert \ll p$, one can essentially ignore
the rows outside $\Omega^*$ using an argument similar to that in the
proof of Theorem \ref{theoremblockdesign} and the problem becomes
equivalent to $\mathrm{ANOVA}(r_*)$ designs. However, condition
\textup{(C3)} allows for the cases $\llvert \bigcup_{i \notin\Omega
^*}S_i\rrvert \gg p$. For example, if $Q=\log(p)^{b}$ for some $b>0$, then
as long as $r^* \gamma_p \gg p a_p \log(p)^{b}$ for some sequence
$a_p \rightarrow\infty$, one can potentially have $n_* Q_p \gg p$,
and hence the simple reduction of the problem as in proof of
Theorem \ref{theoremblockdesign} is no longer possible. In order to
show that the detection problem still behaves similar to an orthogonal
design, one needs much subtler analysis to ignore the information about
the alternative coming from the subjects corresponding to $\mathbf{G}$
part of the design $\bX$. 
Therefore, condition~\textup{(C3)} allows for a rich class of
correlation structures in $\bX$.


%
%
\begin{table}
\tabcolsep=0pt
\caption{Characteristics of the genotype matrix of uncommon/rare
variants of the Dallas Heart study using the parameters defined in
Definition \protect\ref{wa}}\label{table3variantsremoved}
\begin{tabular*}{\tablewidth}{@{\extracolsep{\fill}}@{}lcccccc@{}}
\hline
\textbf{Demography} & $\bolds{r^{*}}$ & $\bolds{n_{*}}$& $\bolds{Q}$ & $\bolds{\frac{n_* Q^2/r^{*} }{p^{1/4}}}$
                    & $\bolds{\frac{n_* Q^2/r^{*} }{\sqrt{p}}}$ & $\bolds{\frac{n_* Q^2/r^{*} }{\log (p)}}$
\\
\hline
Overall & 148.00 & 25.00 & 2.00 & 0.22 & 0.07 & 0.15 \\
White & \phantom{0}14.00 & \phantom{0}2.00 & 2.00 & 0.18 & 0.06 & 0.13 \\
Black & 142.00 & 19.00 & 2.00 & 0.17 & 0.06 & 0.12 \\
Hispanic & \phantom{0}26.00 & \phantom{0}4.00 & 2.00 & 0.20 & 0.06 & 0.14 \\
\hline
\end{tabular*}
\end{table}

 The genotype matrix of the Dallas Heart study data shown in
Figure~\ref{DHS} provides empirical evidence that the assumptions in
Definition \ref{def1} are reasonable for design matrices in
sequencing data.
Specifically, Table~\ref{table3variantsremoved} provides the values
of the parameters used in Definition \ref{def1} that
were calculated using the
Dallas Heart study data for different subpopulations of the study to
motivate our conditions. Here, we assumed a dominant coding of the
alleles for the rare variants (MAF $< 5\%$). In most cases, whenever a
subject has more than one mutation, it does not have more than $2$
mutations, which effectively yields $Q=2$ in our conditions. The last
three columns of Table~\ref{table3variantsremoved} refer to
condition \textup{(C3)}. In particular, small values in these columns
suggest that the size of $\mathbf{G}$ is much smaller than the
orthogonal part of the design, supporting condition \textup{(C3)}. 

In subsequent sections, we study the role of the parameter vector
($n^*$, $n_*$, $r^*$, $r_*$, $Q_{n,p}$, $\gamma_{n,p}$) in deciding the detection
boundary. We first present the analysis of relatively simpler $\mathrm
{ANOVA}$ designs followed by the study of Weakly Correlated Designs.
The analysis of simpler $\mathrm{ANOVA}$ designs provides the crux of
insight for the study of detection boundary under Weakly Correlated
Designs, and at the same time yields cleaner results for easier
interpretation. We will demonstrate that the quantity $\frac{1}{r}$ is
the design\vspace*{1pt} matrix sparsity index when $\bX\in\operatorname
{ANOVA}(r)$. In the case of
Weakly Correlated Designs, $r^*$ and $r_*$ play the same role as that
of $r$ in $\operatorname{ANOVA}(r)$ designs. We divide our study of
each design into two
main sections, namely the Dense Regime $(\alpha\leq\frac{1}{2})$ and
the Sparse Regime $(\alpha> \frac{1}{2})$. In the next section, we
first introduce the tests which will be essential for attaining the
optimal detection boundaries in dense and sparse regimes, respectively.

\section{Tests}\label{sec5}\label{sectests}
We propose in this section the Generalized Likelihood Ratio Test and a
Higher Criticism Test for binary regression models.
We begin by defining \mbox{$Z$-}statistics for Weakly Correlated Designs which
will be required for introducing and analyzing upper bounds later. 
Also, in order to separate the information about\vadjust{\goodbreak} the alternative coming
from the $\mathbf{G}$ part of $\bX$, we define a $Z$-statistic
separately for the nonorthogonal part. With this in mind, we have the
following definitions.
%

\begin{defn}
\label{definitionzstatistic}
Let $\bX$ be a Weakly Correlated Design as in Definition~\ref{def1}.
%
\begin{longlist}[1.]
\item[1.] Define the $j$th $Z$-statistic as follows:
\[
Z_j= \sum_{i \in\Omega_j^*}y_i,\qquad j=1,
\ldots,p.
\]

\item[2.] Letting $\mathbf{G}=\{\mathbf{G}_{ij}\}_{n_*\times p}$ define
\[
Z_j^{\mathbf{G}}=\sum_{i=n-n_*+1}^n
\mathbf{G}_{ij} y_i,\qquad j=1,\ldots,p.
\]
\end{longlist}
\end{defn}


With these definitions, we are now ready to construct our tests.

\subsection{The Generalized Likelihood Ratio Test \texorpdfstring{($\mathrm{GLRT}$)}{(GLRT)}}\label{sec5.1}
We now introduce a test that will be used to attain the detection
boundary in the dense regime. Let $Z_j$ be the $j$th
$Z$-statistic in Definition \ref{definitionzstatistic}. Then the
Generalized Likelihood Ratio Test is based on the following test statistic:
%
%
\begin{equation}
\label{eqnglrttest} \tglrtnew:= \sum_{j=1}^p {
\frac{4(Z_j-({r_j}/{2}))^2}{r_j}}.
\end{equation}

Under $H_0$, we have $\mathbb{E}_{H_0}(\tglrtnew)=p$ and $\operatorname
{Var}_{H_0}(\tglrt)=O(p)$. Hence, $\frac{\tglrt-p}{\sqrt{2p}}$ is
tight. Our test rejects the null when
\[
\frac{\tglrt-p}{\sqrt{2p}}>t_p
\]
for
a suitable $t_p$ to be decided later.

Note that this test only uses partial information from the data. Since
we shall show that, asymptotically using this partial information is
sufficient, we will not lose power in an asymptotic sense. However,
from finite sample performance point of view, it is more desirable to
use the following test using all the data by incorporating information
from $\mathbf{G}$ as well. This test can be viewed as a combination of
$\mathrm{GLRT}$ statistics using the orthogonal and nonorthogonal
parts of $\bX$, respectively. Specifically, we reject the null
hypothesis
\[
\mbox{when:}\qquad \max \biggl\{\frac{\tglrt-p}{\sqrt{2p}},\frac{\sum_{j=1}^p{
[(Z_j^{\mathbf{G}})^2-\mathbb
{E}_{H_0}((Z_j^{\mathbf{G}})^2) ]}}{\sqrt{\mathrm
{V}_{H_0}(\sum_{j=1}^p(Z_j^{\mathbf{G}})^2)}} \biggr
\}>t_p.
\]

Note that given a particular $\mathbf{G}$, the quantities $\mathbb
{E}_{H_0}\{(Z_j^{\mathbf{G}})^2\}$ and\break $\mathrm{V}_{H_0}\{\sum_{j=1}^p(Z_j^{\mathbf{G}})^2\}$ can be easily calculated by simple
moment calculations of Bernoulli random variables. We do not go into
specific details here. Finally, since combining correct size tests by
Bonferroni correction does not change asymptotic power, our proofs
about asymptotic power continue to hold for this modified $\mathrm
{GLRT}$ without any change.

\subsection{Extended Higher Criticism Test}\label{sec5.2} Assume $r_*
\geq2$.
Let $R_j$ be a generic $\operatorname{Bin}(r_j,\frac{1}{2})$ random
variable and $\mathbb{B}_j,\overline{\mathbb{B}_j}$, respectively,
denote the distribution function and the survival function of $\frac
{\llvert R_j-({r_j}/{2})\rrvert }{\sqrt{{r_j}/{4}}}$. Hence,
\[
\mathbb{B}_j(t)=\mathbb{P} \biggl(\frac{\llvert R_j-({r_j}/{2})\rrvert
}{\sqrt{{r_j}/{4}}}\leq t \biggr),
\qquad \overline{\mathbb{B}_j}(t)=1-\mathbb {B}_j(t) .
\]
%
From Definition \ref{definitionzstatistic}, the $Z_j$'s are
independent $\operatorname{Bin}(r_j,\frac{1}{2})$ under $H_0$ for $j
=1,\ldots,p$. 
Let
\[
W_p(t)=\frac{\sum_{j=1}^p{\I({\llvert Z_j-
({r_j}/{2})\rrvert }/{\sqrt{{r_j}/{4}}}> t )-\overline{\mathbb
{B}_j}(t)}}{\sqrt{\sum_{j=1}^p\overline{\mathbb
{B}_j}(t)(1-\overline{\mathbb{B}_j}(t))}}.
\]

Now we 
define the Higher Criticism Test as
%
%
\begin{equation}
\label{eqnthcwa} \thcnew:=\max_{t \in[1,\sqrt{3\log(p)}] \cap\mathbb{N}}W_p(t),
\end{equation}
where $\mathbb{N}$ denotes the set of natural numbers. The next
theorem provides the rejection region for the Higher Criticism Test.
%

\begin{theorem}\label{theoremhighercriticismnullwa}
For Weakly Correlated Designs, $\lim_{p \rightarrow\infty} \mathbb
{P}_{H_0}(\thc>\break  \log(p)) = 0$.
\end{theorem}

Hence, one can use ${(1+\varepsilon)\log(p)}$ as a cutoff to construct a
test based on $\thc$ for any arbitrary fixed $\varepsilon>0$:
%
%
\begin{equation}
\label{highercriticismtestwa}
\mbox{Higher Criticism Test:\qquad Reject when } \thcnew>
{(1+\varepsilon)\log(p)}.
\end{equation}
%
By Theorem \ref{theoremhighercriticismnullwa}, the above test
based on $\thcnew$ has asymptotic type I error converging to $0$.
We note that, when $r_*\gg\log(p)$, we can obtain a rejection region
of the form $\thcnew> \sqrt{{2(1+\varepsilon)\log{\log(p)}}}$ while
maintaining asymptotic type I error control. This type of rejection
region is common in the Higher Criticism literature. As we will see in
Section~\ref{secSADE}, the interesting regime where the Higher
Criticism Test is important is when $r_*\gg\log(p)$. In this regime,
we can have the same rejection region of the Higher Criticism as
obtained in \citet{Jin1} and \citet{Jin2}. However, for generality we will
instead work with the rejection region given by equation (\ref
{highercriticismtestwa}).

Note that this test only uses partial information from the data. We
shall show that, asymptotically, using this partial information is
sufficient, we will not lose power in an asymptotic sense. However,
from a finite sample performance point of view, it is more desirable to
use the following test using all the data by incorporating information
from $\mathbf{G}$. The below can be viewed as a combination of Higher
Criticism Tests based on the orthogonal and nonorthogonal parts of $\bX
$, respectively.

Specifically, letting $g_j=\sum_{i>n^*}\bX_{ij}, j=1,\ldots,p$,
define the Higher Criticism type test statistic based on $\mathbf{G}$ as
\begin{eqnarray*}
&& W^{\mathbf{G}}_p(t)
\\
&&\qquad =\frac{\sum_{j=1}^p{\I(
{\llvert Z_j^{\mathbf{G}}-({g_j}/{2})\rrvert }/{\sqrt{{g_j}/{4}}}> t
)-\mathbb{P}_{H_0} ({\llvert Z_j^{\mathbf{G}}-({g_j}/{2})\rrvert }/{\sqrt{{g_j}/{4}}}> t )}}{\sqrt{\operatorname
{Var}_{H_0}\sum_{j=1}^p{\I({\llvert Z_j^{\mathbf
{G}}-({g_j}/{2})\rrvert }/{\sqrt{{g_j}/{4}}}> t )}}}.
\end{eqnarray*}
The quantities $\mathbb{P}_{H_0} \{\frac{\llvert Z_j^{\mathbf{G}}-({g_j}/{2})\rrvert }{\sqrt{{g_j}/{4}}}> t \}$ and $\operatorname
{Var}_{H_0}\sum_{j=1}^p{\I\{\frac{\llvert Z_j^{\mathbf
{G}}-({g_j}/{2})\rrvert }{\sqrt{{g_j}/{4}}}> t \} }$ can be
suitably approximated based on $\mathbf{G}$. However, we omit the
specific details here for coherence of exposition. Finally, defining
\[
W_p^{\mathrm{comb}}(t)=\max \bigl\{ {W_p(t),W^{\mathbf{G}}_p(t)}
\bigr\},
\]
one can follow the previous steps in defining the Higher Criticism Test
with exactly similar arguments. Since combining correct size tests by
Bonferroni correction does not change asymptotic power, the proofs
concerning the power of the resulting test goes through with similar
arguments. We omit the details here.

\section{Detection boundary and asymptotic analysis for \texorpdfstring{$\mathrm{ANOVA}$}{ANOVA} designs}\label{sec6} \label{secSADE}
We begin by noting that the
$\operatorname{ANOVA}(r)$
designs can be equivalently cast as a problem of testing homogeneity
among $p$ different binomial populations with $r$ trials each. Suppose
%
%
\begin{equation}
\label{eqnbinomialproportionmodel} y_j \sim\operatorname{Bin} \bigl(r,
\tfrac{1}{2}+ \nu_j \bigr) \qquad\mbox{independent for } j=1,
\ldots,p.
\end{equation}
Let $\bnu=(\nu_1,\ldots,\nu_p)^t$. For some $ \Delta\in(0,\frac
{1}{2}]$, we are interested in testing the global null hypothesis
%
%
\begin{equation}
\label{eqnbinomialproportionhypothesis} H_0\dvtx \bnu=0\quad\mbox {vs}\quad H_1
\dvtx \bnu\in \Xi^{\Delta}_k= \bigl\{ \bnu\in
R_k^p\dvtx \min \bigl\{\llvert\nu_j\rrvert
\dvtx \nu_j \neq0 \bigr\} \geq\Delta \bigr\}.
\end{equation}
When $\bX\in\operatorname{ANOVA}(r)$, models \eqref{eqnmodmain}
and \eqref
{eqnbinomialproportionmodel} are equivalent with $\eta_j=\theta
(\beta_j)-\frac{1}{2}$. Hence, sparsity in $\bbeta$ is equivalent to
sparsity in $\bnu$ in the sense that $\bbeta\in R_k^p$ if and only if
$\bnu\in R_k^p$. Further, the rate of $\Delta$, which determines the
asymptotic detectability of \eqref
{eqnbinomialproportionhypothesis}, can be related to the rate of
$A$, which determines detectability in \eqref{eqnhyp} when the link
function $\theta$ is continuously differentiable in a neighborhood
around $0$.

%
\begin{rem}
\label{remarkrelationdetectionboundary}
When $\theta$ is the distribution function for a uniform random
variable $\mathrm{U}(-\frac{1}{2},\frac{1}{2})$, then $\nu_j=\beta
_j$ for all $j=1,\ldots,p$. Hence, the detection boundary in problem
\eqref{eqnbinomialproportionhypothesis} follows from that in
problem \eqref{eqnhyp} by taking\vspace*{1pt} $\theta$ to be the distribution
function of $\mathrm{U}(-\frac{1}{2},\frac{1}{2})$, that is,
$\theta(x)=(x+\frac{1}{2}) \I(-\frac{1}{2}<x<\frac{1}{2})$.
\end{rem}

%
\begin{rem}
The prior $\pi_{\mathrm{eq}}$ that we will use for testing for the binomial
homogeneity of proportions is as follows. For each $m \in M(k, p)$, let
$\xi^m =(\xi_j)_{j \in m}$ be a sequence of independent Rademacher
random variables taking values in $\{+1,-1\}$ with equal probability.
Given $\Delta\in(0,\frac{1}{2})$ for testing \eqref
{eqnbinomialproportionhypothesis}, a realization from the prior
distribution $\pi_{\mathrm{eq}}$ on $\mathbb{R}^p$ can be expressed
as $\nu_{\xi,m} =\sum_{j \in m}{\Delta\xi_j e_j}$, where
$(e_j)_{j=1}^p$ is the canonical basis of $\mathbb{R}^p$ and $m$ is
uniformly chosen from $M(k,p)$. Note that given the prior $\pi$ on
$\bbeta=(\beta_1,\ldots,\beta_p)^r$ discussed earlier, $\pi
_{\mathrm{eq}}$ is the prior induced on $\bnu=(\nu_1,\ldots,\nu
_p)^t$ with $\frac{1}{2}+\nu_j=\theta(\beta_j)$ for $j=1,\ldots,p$.
\end{rem}

Owing to Remark \ref{remarkrelationdetectionboundary}, one can
deduce the detection boundary of the binomial proportion model \eqref
{eqnbinomialproportionmodel} from the detection boundary in
$\operatorname{ANOVA}(r)$
designs. However, for the sake of easy reference, we provide the
detection boundaries for both models. Before proceeding further, we
first state a simple result about $\mathrm{ANOVA}$ designs, a part of
which directly follows from Theorem \ref{theoremnodetectcondition}.
Note that $\operatorname{ANOVA}(1)$ design corresponds to the case when the
design matrix is identity $I_{p \times p}$. Unlike Gaussian linear
models, for binary regression, when the design matrix is identity, for
two-sided alternatives, all tests are asymptotically powerless
irrespective of sparsity (i.e., in both dense and sparse
regimes) and signal strengths. Such a result arises for $r=1$ because
we allow the alternative to be two-sided.
In the modified problem where one only considers the one-sided
alternatives, all tests still remain asymptotically powerless
irrespective of signal strengths when $r=1$ in the sparse regime,
{that is}, when $\alpha>\frac{1}{2}$. However, in the dense
regime, {that is}, when $\alpha\leq\frac{1}{2}$, the problem
becomes nontrivial and the test based on the total number of successes
attains the detection boundary. The detection boundary for this
particular problem is provided in Theorem \ref{identitydesign} part
2(b). Also, in the one-sided problem, the Bayes test can be explicitly
evaluated and quite intuitively turns out to be a function of the total
number of successes. In the next theorem, we collect all these results.

\begin{theorem}\label{identitydesign}
Assume $\bX\in\operatorname{ANOVA}(1)$, which assumes $r=1$ and $\bX=I$.
Then the following holds for both problems \eqref{eqnhyp} and \eqref
{eqnbinomialproportionhypothesis}.
\begin{enumerate}
\item For two-sided alternatives all tests are asymptotically powerless
irrespective of sparsity and signal strength.
\item For one-sided alternatives:
\begin{longlist}[(a)]
\item[(a)] Suppose $\theta\in\BC^1(0)$, which is defined in Section~\ref{sec2.1}.
Then in the dense regime ($\alpha\leq\frac{1}{2}$), all tests are
asymptotically powerless if $\frac{A^2}{p^{1-2\alpha}} \rightarrow0$
in problem~\eqref{eqnhyp} or $\frac{\Delta^2}{p^{1-2\alpha}}
\rightarrow0$ in problem~\eqref{eqnbinomialproportionhypothesis}.
Further, if $\frac{A^2}{p^{1-2\alpha}} \rightarrow\infty$ in
problem~\eqref{eqnhyp} or $\frac{\Delta^2}{p^{1-2\alpha}}
\rightarrow\infty$ in problem~\eqref
{eqnbinomialproportionhypothesis}, then the test based on the total
\mbox{number} of successes ($\sum_{i=1}^p y_i$) is asymptotically powerful.
\item[(b)] In sparse regime ($\alpha> \frac{1}{2}$), all tests are
asymptotically powerless.
\end{longlist}
\end{enumerate}
\end{theorem}

The case of two-sided of alternatives when $r=1$ can indeed be
understood in the following way. Under the null hypothesis, each $y_i$
is an independent Bernoulli$(1/2)$ random variable and under the prior
on the alternative which allows each $\beta_i$ to be $+A$ or $-A$ with
probability $\frac{1}{2}$, the $y_i$'s are again independent
Bernoulli$(1/2)$ random variables. So, of course, there is no way to
distinguish them based on the observations $y_i$'s when the $\bbeta$
is generated according to the prior mentioned earlier. Our proof is
based on this heuristic. However, the above argument is invalid even
for $r>1$ and one can expect nontrivial detectability conditions on $A$
when $r>1$. In the dense regime, we observe that simply $r>1$ is enough
for this purpose. However, the sparse regime requires a more delicate
approach in terms of the effect of $r>1$.
%

\begin{rem}
Note that Theorem \ref{identitydesign}, other than part 2(a),
requires no additional assumption on $\theta$ other than the symmetry
requirement in equation \eqref{eqntheta}.
\end{rem}

\subsection{Dense regime \texorpdfstring{$(\alpha\leq\frac{1}{2})$}{(alphaleq{1}/{2})}}\label{sec6.1}
The detection complexity in the dense regime with $r> 1$ matches the
Gaussian linear model case. Interestingly, just by increasing $1$
observation per treatment from the identity design matrix scenario, the
detection boundary changes completely. 
The following theorem provides the lower and upper bound for the dense
regime when $r>1$.
%

\begin{theorem}\label{theorembalancedlargerdense} Let $\bX\in
\operatorname{ANOVA}(r)
$. Let $k=p^{1-\alpha}$ with $\alpha\leq\frac{1}{2}$ and the block
size/binomial denominator $r > 1$.
\begin{enumerate}[2.]
\item Consider the model \eqref{eqnmodmain} and the testing problem
given by \eqref{eqnhyp}. Assume $\theta\in\BC^1(0)$. Then:
\begin{longlist}[(a)]
\item[(a)] If $A \ll\sqrt{\frac{p^{1/2}}{kr}}$, then all tests are
asymptotically powerless.\vspace*{1pt}

\item[(b)] If $A \gg\sqrt{\frac{p^{1/2}}{kr}}$, then the $\mathrm{GLRT}$
is asymptotically powerful.
\end{longlist}
\item Consider model \eqref{eqnbinomialproportionmodel} and the
testing problem \eqref{eqnbinomialproportionhypothesis}. Then:
\begin{longlist}[(a)]
\item[(a)] If $\Delta\ll\sqrt{\frac{p^{1/2}}{kr}}$, then all tests are
asymptotically powerless.\vspace*{1pt}
\item[(b)] If $\Delta\gg\sqrt{\frac{p^{1/2}}{kr}}$, then the $\mathrm
{GLRT}$ is asymptotically powerful.
\end{longlist}
\end{enumerate}
\end{theorem}

Also\vspace*{2pt} when $\frac{A^2 kr}{\sqrt{p}}$ or $\frac{\Delta^2 kr}{\sqrt
{p}}$ remains bounded away from $0$ and $\infty$, the asymptotic power
of $\mathrm{GLRT}$ remains bounded between $0$ and $1$. The upper and
lower bound rates of the minimum signal strength match with that of
\citet{Candes} and \citet{Ingster5}.

\subsection{Sparse regime \texorpdfstring{$(\alpha>\frac{1}{2})$}{(alpha>{1}/{2})}}\label{sec6.2}
Unlike the dense regime, the sparse regime depends more heavily on the
value of $r$. 
The next theorem quantifies this result; it shows that in the sparse
regime if $r \ll\log(p)$, then all tests are asymptotically
powerless. Indeed this can be argued from Theorems \ref
{theoremnodetectcondition} and \ref{theoremblockdesign}. However,
for the sake of completeness, we provide it here.

%
\begin{theorem}\label{theorembalancedsmallr}
Let $k=p^{1-\alpha}$ with $\alpha> \frac{1}{2}$. If $r \ll\log
(p)$, then for both the problems and \eqref{eqnhyp} and \eqref
{eqnbinomialproportionhypothesis}, all tests are asymptotically powerless.
\end{theorem}

%
\begin{rem}
Theorem \ref{theorembalancedsmallr} requires no additional
smoothness assumption on~$\theta$ other than the symmetry requirement
in equation \eqref{eqntheta}.
\end{rem}

Thus, for the rest of this section we consider the case where $k \ll
\sqrt{p}$ and $r \gg\log(p)$. We first divide our analysis into two
parts, where we study the lower bound and upper bound of the problem separately.
\subsubsection{Lower bound}\label{sec6.2.1}
\begin{figure}

\includegraphics{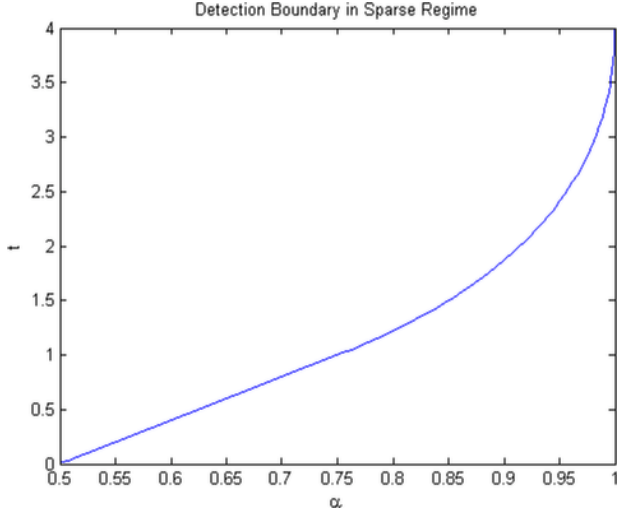}

\caption{Detection boundary $t =\rho_{\mathrm{binary}}^{*}(\alpha)$
in the sparse regime when $\theta$ corresponds to logistic regression.
The detectable region is $t >\rho_{\mathrm{binary}}^{*}(\alpha)$,
and the undetectable region is $t <\rho_{\mathrm{binary}}^{*}(\alpha
)$. The curve corresponds to $t=\rho_{\mathrm{binary}}^{*}(\alpha
)$.} \label{fig1}
\end{figure}

To introduce a sharp lower bound in the regime where $\alpha>\frac
{1}{2}$ and $r \gg\log(p)$ in the binary regression model \eqref
{eqnmodmain} and the testing problem \eqref{eqnhyp} for the
ANOVA($r$) design, we define the following functions. Figure~\ref{fig1}
provides a graphical representation of the detection boundary. Define
%
%
\begin{equation}
\rho_{\mathrm{binary}}^{*}(\alpha) = \cases{ \displaystyle
\frac{(\alpha-({1}/{2}))}{4(\theta'(0))^2}, &\quad if $\displaystyle \frac{1}{2} <\alpha<
\frac{3}{4}$,
\vspace*{5pt}\cr
\displaystyle\frac{(1-\sqrt{1-\alpha})^2}{4(\theta
'(0))^2},&\quad if $\displaystyle\alpha\geq\frac{3}{4}$.}
\end{equation}

%

This\vspace*{1pt} is the same as the Gaussian detection boundary \eqref
{eqndetboundLM} multiplied by $1/4(\theta'(0))^2$. The reason
for the appearance of the factor $1/4(\theta'(0))^2$ is that the
Fisher information for a single Bernoulli\vspace*{1pt} sample under binary
regression model~\eqref{eqnmodmain} is equal to $\sqrt{4(\theta
'(0))^2}$.

For every $j \in\{1,\ldots,p\}$, we have
\[
\hat{\beta}_j^{\mathrm{MLE}} \stackrel{d} {\rightarrow} N \bigl(
\beta_j,\sigma_j^2 \bigr),
\]
where\vspace*{1pt} $\sigma_j^2=4(\theta'(0))^2$ under $H_0$ and $\sigma^2_j
\approx4(\theta'(0))^2$ under $H_1$. To see this, note that under
$H_1$ we have $\sigma^2_j =(\frac{1}{2}+\delta)(\frac{1}{2}-\delta
)\approx4\theta'(0)$ where $\delta>0$ is small and denotes a
departure of the Bernoulli proportion from the null value\vspace*{1pt} of $\frac
{1}{2}$, {that is}, under $H_1$, the outcomes corresponding to the
signals follow $\operatorname{Bernoulli}(\frac{1}{2}+\delta)$ or
$\operatorname{Bernoulli}(\frac{1}{2}-\delta)$. 
This implies $\sqrt{\frac{1}{4(\theta'(0))^2}} \hat{\bbeta}$
should yield\vspace*{1pt} a detection boundary similar to the multivariate Gaussian
model case.

For the detection boundary in the corresponding binomial proportion
model~\eqref{eqnbinomialproportionmodel} and the testing problem~\eqref{eqnbinomialproportionhypothesis}, we define the following function:
%
%
\begin{equation}
\rho_{\mathrm{binomial}}^{*}(\alpha) = \cases{ \displaystyle
\frac{(\alpha-({1}/{2}))}{4}, &\quad if $\displaystyle\frac{1}{2} <\alpha<
\frac{3}{4}$,
\cr
\displaystyle\frac{(1-\sqrt{1-\alpha})^2}{4}, &\quad if $\displaystyle
\alpha\geq\frac{3}{4}$.}
\end{equation}
%

The following theorem provides the exact lower boundary for the\break
$\operatorname{ANOVA}(r)$
designs for the binary regression model as well as the corresponding
binomial problem.
%

\begin{theorem}\label{theorembalancedlargersparse}
Let $\bX\in\operatorname{ANOVA}(r)$. Suppose $r \gg\log(p)$ and
$k=p^{1-\alpha}$ with
$\alpha> \frac{1}{2}$.
\begin{enumerate}[2.]
\item Consider the binary regression model \eqref{eqnmodmain} and the
testing problem \eqref{eqnhyp}. Further suppose that $\theta\in\BC
^2(0)$. Let $A=\sqrt{\frac{2t\log(p)}{r}}$. If $t<{\rho_{\mathrm
{binary}}^{*}(\alpha)}$, all tests are asymptotically powerless.
\item Consider the binomial model \eqref
{eqnbinomialproportionmodel} and the testing problem \eqref
{eqnbinomialproportionhypothesis}.
Let $\Delta=\sqrt{\frac{2t\log(p)}{r}}$. If $t<{{\rho_{\mathrm
{binomial}}^{*}(\alpha)}}$, all tests are asymptotically powerless.
\end{enumerate}
\end{theorem}

%
\begin{rem}
As mentioned in the \hyperref[sec1]{Introduction}, the analysis turns out to be
surprisingly nontrivial since it seems not possible to simply reduce
the calculations to the Gaussian case by doing a Taylor expansion of
$L_{\pi}$ around $\bbeta=0$. In particular, a natural approach to
analyze these problems is to expand the integrand of $L_{\pi}$ by a
Taylor series around $\bbeta=0$ and thereby reducing the analysis to
calculations in the Gaussian situation and a subsequent analysis of the
remainder term. However, in order to find the sharp detection boundary,
the analysis of the remainder term turns out to be very complicated and
nontrivial. Thus, our proof to Theorem~\ref
{theorembalancedlargersparse} is not a simple application of
results from the Gaussian linear model.
\end{rem}
%
\subsubsection{Upper bound}\label{sec6.2.2}
According to Theorem \ref{theorembalancedlargersparse}, all tests
are asymptotically powerless if $t<\rho_{\mathrm{binary}}^{*}(\alpha
)$ in the sparse regime. In this section, we introduce tests which
reach the lower bound discussed in the previous section. We divide our
analysis into two subsections.
In Section \ref{sec6.2.2.1}, we study the Higher Criticism Test defined by
\eqref{eqnthcwa} which is asymptotically powerful as soon as $t>\rho
_{\mathrm{binary}}^{*}(\alpha)$ in the sparse regime. In Section
\ref{sec6.2.2.3}, we discuss a more familiar Max Test or minimum $p$-value test
which attains the sharp detection boundary only for $\alpha\geq\frac
{3}{4}$. 

\paragraph{The Higher Criticism Test}\label{sec6.2.2.1}
In this section, we study the version of Higher Criticism introduced in
Section \ref{sec6.2}. Recall, 
we have by Theorem \ref{theoremhighercriticismnullwa} that the
type I error of the Higher Criticism Test, as defined by equation (\ref
{highercriticismtestwa}), converges to $0$.
The next theorem states the optimality of the Higher Criticism Test as
soon as the signal strength exceeds the detection boundary.

%
\begin{theorem}\label{theoremhighercriticism}
Let $\bX\in\operatorname{ANOVA}(r)$. Suppose $r \gg\log(p)$ and
$k=p^{1-\alpha}$ with
$\alpha> \frac{1}{2}$.
\begin{enumerate}[2.]
\item Consider the binary regression model \eqref{eqnmodmain} and
the testing problem \eqref{eqnhyp}. Further suppose that $\theta\in
\BC^2(0)$. Let $A=\sqrt{\frac{2t\log(p)}{r}}$. If\vspace*{1pt} $t>{\rho
_{\mathrm{binary}}^{*}(\alpha)}$, then the Higher Criticism Test is
asymptotically powerful.

\item Consider the binomial model \eqref
{eqnbinomialproportionmodel} and the testing problem \eqref
{eqnbinomialproportionhypothesis}.
Let $\Delta=\sqrt{\frac{2t\log(p)}{r}}$. If $t>{{\rho_{\mathrm
{binomial}}^{*}(\alpha)}}$, then the Higher Criticism Test is
asymptotically powerful.
\end{enumerate}
\end{theorem}


\paragraph{Comparison with the original Higher Criticism Test}\label{sec6.2.2.2}
We begin by providing a slight simplification of $\thc$ in
$\operatorname{ANOVA}(r)$
designs. Let $S$ be a generic $\operatorname{Bin}(r,\frac{1}{2})$ random
variable and $\mathbb{B},\overline{\mathbb{B}}$, respectively,
denote the distribution function and the survival function of $\frac
{\llvert S-({r}/{2})\rrvert }{\sqrt{{r}/{4}}}$. Hence,
\[
\mathbb{B}(t)=\mathbb{P}\biggl(\frac{\llvert S-({r}/{2})\rrvert }{\sqrt
{{r}/{4}}}\leq t\biggr),\qquad \overline{
\mathbb{B}}(t)=1-\mathbb{B}(t) .
\]
In\vspace*{1pt} the case of $\operatorname{ANOVA}(r)$ designs, $W_p(t)=\sqrt
{p}\frac{\overline
{\mathbb{F}}_p(t)-\overline{\mathbb{B}}(t)}{\sqrt{\overline
{\mathbb{B}}(t)(1-\overline{\mathbb{B}}(t))}}$.
The original Higher Criticism Test as defined by \citet{Jin1} can also
be calculated as a maximum over some appropriate function of $p$-values.
By that token, ideally we would like to define the Higher Criticism
Test statistic as
\[
\thc^{\mathrm{Ideal}}=\sup_{0<t<{r}/{2}}{W_p(t)}.
\]
%

However, due to difficulties in calculating the null distribution for
deciding a cut-off for the rejection region, we instead work with a
discretized version of it. We detail this below in the context of
$\operatorname{ANOVA}(r)
$ designs. Define the $j$th $p$-value as $q_j=\mathbb
{P}(\llvert \operatorname{Bin}(r,\frac{1}{2})-\frac{r}{2}\rrvert
>\llvert Z_j-\frac{r}{2}\rrvert )$
for $1,\ldots,p$ and let $q_{(1)},\ldots,q_{(p)}$ be the ordered
$p$-values based on exact binomial distribution probabilities. Define
\[
\thc' =\max_{1\leq j \leq p}\sqrt{p}\frac{({j}/{p})-q_{(j)}}{\sqrt{q_{(j)}(1-q_{(j)})}}.
\]
It is difficult to analyze the distribution of $\thc'$ under the null
to decide a valid cut-off for testing. The following proposition yields
a relationship between $\thc$, $\thc^{\mathrm{Ideal}}$ and $\thc'$.
%

\begin{prop}\label{lemhcor}
Let\vspace*{1pt} $\llvert Z-\frac{r}{2}\rrvert _{(j)}$ denote the $j$th order
statistics based on $\llvert Z_i-\frac{r}{2}\rrvert, i=1,\ldots,p$.
For $t$ such
that $\llvert Z-\frac{r}{2}\rrvert _{(p-j)}\leq t<\llvert Z-\frac
{r}{2}\rrvert _{(p-j+1)}$, we have
\[
\sqrt{p}\frac{\overline{\mathbb{F}}_p(t)-\overline{\mathbb
{B}}(t)}{\sqrt{\overline{\mathbb{B}}(t)(1-\overline{\mathbb
{B}}(t))}} \leq\sqrt{p}\frac{({j}/{p})-q_{(j)}}{\sqrt
{q_{(j)}(1-q_{(j)})}}.
\]
\end{prop}

Hence, from Proposition~\ref{lemhcor}, we observe that we have the
following inequality:
%
%
\begin{equation}
\label{eqnidealhighercriticismcomparisions} \thc' \geq\thc^\mathrm {Ideal} \geq\thc.
\end{equation}
This unlike the results in \citet{Jin1} and \citet{Cai1}, where the
leftmost inequality is a equality. Therefore, it is worth further
comparing the above discussion to the Higher Criticism Test introduced
by \citet{Jin1,Jin2} in the Gaussian framework. In the case of
orthogonal Gaussian linear models, $\thc, \thc'$ and $\thc^\mathrm
{Ideal}$ are defined by standard normal survival functions and
$p$-values, respectively, and one uses $Z_j$ instead\vspace*{1pt} of $\frac{Z_j-({r}/{2})}{\sqrt{{r}/{4}}}$ in the definition of $\thc$. This
yields that in the Gaussian framework the leftmost inequality of \eqref
{eqnidealhighercriticismcomparisions} is an equality. Moreover,
under the framework, standard empirical process results for continuous
distribution functions yield asymptotics for $\thc^\mathrm{Ideal}$
under the null. Therefore, in the Gaussian case the uncountable
supremum in the definition of $\thc^\mathrm{Ideal}$ is attained and
the statistic is algebraically equal to a maximum over finitely many
functions of $p$-values, namely, $\thc'$. However, due to the
possibility of strict inequality in Proposition \ref{lemhcor} for the
binomial distribution, we cannot reduce our computation to $p$-values as
in the Gaussian case. Although it is true that marginally each $q_j$ is
stochastically smaller than a $\mathrm{U}(0,1)$ random variable, we
are unable to find a suitable upper bound for the rate of $\thc'$
since it also depends on the joint distribution of $q_{(1)},\ldots
,q_{(p)}$. It might be possible to estimate the gaps between $\thc
', \thc^\mathrm{Ideal}$ and $\thc$, but since this is not
essential for our purpose, we do not attempt this. 


\paragraph{Rate optimal upper bound: Max Test}\label{sec6.2.2.3}
A popular multiple comparison procedure is the minimum $p$-value test. In
the context of Gaussian linear regression, \citet{Jin1} and \citet
{Candes} showed that the minimum $p$-value test reaches the sharp
detection boundary if and only if $\alpha\geq\frac{3}{4}$. In this
section, we introduce and study the minimum $p$-value test in binary
regression models.

As before, define the $j$th $p$-value as
\[
q_j=\mathbb{P} \biggl( \biggl\llvert\operatorname{Bin} \biggl(r,
\frac{1}{2} \biggr)-\frac{r}{2} \biggr\rrvert> \biggl\llvert
Z_j-\frac{r}{2} \biggr\rrvert \biggr)
\]
for $j = 1,\ldots,p$ and let $q_{(1)},\ldots,q_{(p)}$ be the ordered
$p$-values. We will study the test based on the minimum $p$-value
$q_{(1)}$. Note that it is equivalent to study the test based on the statistic
\[
\max_{1 \leq j \leq p}{W_{j}}, \qquad W_{j}=
\frac{\llvert Z_j-({r}/{2})\rrvert }{\sqrt{{r}/{4}}}.
\]
From now on, we will call this the Max Test. In the following theorem,
we show that similar to Gaussian linear models, for binary regression,
the\vspace*{1pt} Max Test attains the sharp detection boundary if and only if
$\alpha\geq\frac{3}{4}$. However, if one is interested in rate
optimal testing, {that is}, only the rate or order of the
detection boundary rather than the exact constants, the Max Test
continues to perform well in the entire sparse regime.
%

\begin{theorem}\label{theoremmaxtest}
Let $\bX\in\operatorname{ANOVA}(r)$. Suppose $r \gg(\log{r})^2
\log(p)$ and
$k=p^{1-\alpha}$ with $\alpha> \frac{1}{2}$.
\begin{enumerate}[2.]
\item Suppose $\theta\in\BC^2(0)$ and let $A=\sqrt{\frac{2t\log
(p)}{r}}$. Set
\[
\rho^{*}_\mathrm{Max,\mathrm{binary}}(\alpha)=
\frac{(1-\sqrt
{1-\alpha})^2}{4(\theta'(0))^2}.
\]
Then in the model \eqref{eqnmodmain} and problem \eqref{eqnhyp} one
has the following:
\begin{longlist}[(a)]
\item[(a)] If $t>\rho^{*}_\mathrm{Max,\mathrm{binary}}(\alpha)$, then
the Max Test is asymptotically powerful.
\item[(b)] If $t<\rho^{*}_\mathrm{Max,\mathrm{binary}}(\alpha)$, then
the Max Test is asymptotically powerless.
\end{longlist}
\item Let $\Delta=\sqrt{\frac{2t\log(p)}{r}}$. Set
$\rho^{*}_\mathrm{Max,\mathrm{binomial}}(\alpha)=\frac{(1-\sqrt
{1-\alpha})^2}{4} $.
Then in the model \eqref{eqnbinomialproportionmodel} and problem
\eqref{eqnbinomialproportionhypothesis} one has the following:
\begin{longlist}[(a)]
\item[(a)] If $t>\rho^{*}_\mathrm{Max,\mathrm{binomial}}(\alpha)$, then
the Max Test is asymptotically powerful.

\item[(b)] If $t<\rho^{*}_\mathrm{Max,\mathrm{binomial}}(\alpha)$, then
the Max Test is asymptotically powerless.
\end{longlist}
\end{enumerate}
\end{theorem}

Theorem \ref{theoremmaxtest} implies that the detection boundary for
the Max Test matches the detection boundary of the Higher Criticism
Test only for $\alpha\geq\frac{3}{4}$. For $\alpha< \frac{3}{4}$,
the detection boundary of the Max Test lies strictly above that of the
Higher Criticism Test. Hence, the Max Test fails to attain the sharp
detection boundary in the moderate sparsity regime, $\alpha< \frac
{3}{4}$. Thus, if one is certain of high sparsity it can be reasonable
to use the Max Test whereas the Higher Criticism Test performs well
throughout the sparse regime. It is worth noting that the requirement
$r \gg(\log(r))^2 \log(p)$ is a technical constraint and can be
relaxed. In most situations, it does not differ much from the actual
necessary condition $r \gg\log(p)$, and hence we use $r \gg(\log
(r))^2 \log(p)$ for proving Theorem \ref{theoremmaxtest}.

\section{Detection boundary and asymptotic analysis for Weakly Correlated Designs}\label{sec7}
\label{secWA}
In this section, we study the role of the parameter vector
$(n^*,n_*,\break  r^*,r_*, Q_{n,p},\gamma_{n,p})$ in deciding the detection
boundary for Weakly Correlated Designs defined in Definition \ref{wa}. For
the sake of brevity, we will often drop the subscripts $n,p$ from $Q$
and $\gamma$ when there is no confusion. Recall $\Omega^{*}$ from
Section~\ref{secDD}.

If we just concentrate on the observations corresponding to the rows
with index in $\Omega^{*}$, we have an orthogonal design matrix
similar to $\operatorname{ANOVA}(r)$ designs. Our proofs of lower
bounds in both dense and
sparse regimes and also the test statistics proposed for the attaining
the sharp upper bound is motivated by this fact. Similar to
$\operatorname{ANOVA}(r)$
designs, we divide our analysis into the dense and sparse regimes.
Also, owing to the possible nonorthogonality of $\bX$ for Weakly
Correlated Designs, we cannot directly reduce this problem to testing
homogeneity of binomial proportions as in \eqref
{eqnbinomialproportionhypothesis}. So, henceforth, we will be
analyzing model \eqref{eqnmodmain} and corresponding testing problem
\eqref{eqnhyp}. However, as we shall see, under certain combinations
of $(n^*,n_*,r^*,r_*,Q,\gamma)$, one can essentially treat the problem
as an orthogonal design like in $\operatorname{ANOVA}(r)$ designs.
This is explained in
the following two sections.

\subsection{Dense regime \texorpdfstring{$(\alpha\leq\frac{1}{2})$}{(alphaleq{1}/{2})}}\label{sec7.1} %
We recall the definition of the GLRT from equation \eqref
{eqnglrttest}. The following theorem provides the lower and upper
bound for the dense regime.
%

\begin{theorem}\label{theoremwalargerdense}
Let $\bX$ be a Weakly Correlated Design as in Definition \ref{def1}.
Suppose Let $k=p^{1-\alpha}$ with $\alpha\leq\frac{1}{2}$ and $r_*
> 1$. Assume $\theta\in\BC^2(0)$ and set $\gamma=p^{(1/2)-\alpha}$. Then we have the following:
\begin{longlist}[2.]
\item[1.] If $A \ll\sqrt{\frac{p^{1/2}}{kr^*}}$, then all tests are
asymptotically powerless.

\item[2.] If $A \gg\sqrt{\frac{p^{1/2}}{kr_*}}$, then the $\mathrm
{GLRT}$ is asymptotically powerful.
\end{longlist}
\end{theorem}

We note that the form of the detection boundary is exactly same as that
in Theorem \ref{theorembalancedlargerdense} for $\operatorname
{ANOVA}(r)$ designs
with $r^*$ and $r_*$ playing the role of $r$. This implies that when
$n_* Q^2$ is not too large ($\frac{n_* Q^2}{r^*} \ll\gamma=p^{(1/2)-\alpha}$); we can still recover the same results as in
$\operatorname{ANOVA}(r)$
designs because the columns of the design matrix are weakly correlated.
\subsection{Sparse regime \texorpdfstring{$(\alpha>\frac{1}{2})$}{(alpha>{1}/{2})}}\label{sec7.2}
Unlike the dense regime, the sparse regime depends more heavily on the
values of $r^*$ and $r_*$. 
The next theorem quantifies this result; it shows that in the sparse
regime if $r^* \ll\log(p)$, then all tests are asymptotically
powerless. This result is analogous to Theorem \ref
{theorembalancedsmallr} for $\operatorname{ANOVA}(r)$ designs.
Indeed this can be
argued from Theorems \ref{theoremnodetectcondition} and \ref
{theoremblockdesign}. However, for the sake of completeness, we
provide it here.
%

\begin{theorem}\label{theoremwasmallr}
Let $\bX$ be a Weakly Correlated Design as in Definition \ref{def1}.
Let $k=p^{1-\alpha}$ with $\alpha> \frac{1}{2}$ and let $\llvert
\bigcup_{i \notin\Omega^*} S_i\rrvert \ll p$. 
If $r^* \ll\log(p)$, then all tests are asymptotically powerless.
\end{theorem}

%
\begin{rem}
The condition $\llvert \bigcup_{i \notin\Omega^*} S_i\rrvert \ll p$,
restricts the location of nonzero elements in the support of rows of
$\bX$ when the row has more than one nonzero element. This restriction
imposes a structure on the deviation of $\bX$ from orthogonality. As
the proof of Theorem \ref{theoremwasmallr} will suggest, this
condition ensures that the assumptions of Theorem \ref
{theoremnodetectcondition} hold, and hence renders all tests
asymptotically powerless irrespective of signal strength. 
\end{rem}

The following theorem provides the value of $\gamma$ that is defined
in condition (C3) in Definition \ref{wa}, to ensure the results parallel to
Theorem \ref{theorembalancedlargersparse}. Not surprisingly, the
test attaining the sharp lower bound turns to be the version of the
Higher Criticism Test introduced in Section~\ref{secSADE}. Similar to
the $\operatorname{ANOVA}(r)$ design, it is also possible to introduce
and study the Max
Test which attains the sharp detection boundary only for $\alpha\geq
\frac{3}{4}$. However, we omit this since it can be easily derived
from the existing arguments.
%

\begin{theorem}\label{theoremwalargersparse}
Let $\bX$ be a Weakly Correlated Design as in Definition \ref{def1}
and $k=p^{1-\alpha}$ with $\alpha> \frac{1}{2}$. Suppose $r_* \gg
\log(p)$, $\gamma=\log(p)$, where $\gamma$ is defined in Definition
\ref{wa}. Further suppose that $\theta\in\BC^2(0)$.
\begin{longlist}[2.]
\item[1.] Let $A=\sqrt{\frac{2t\log(p)}{r^*}}$. If $t<{\rho_{\mathrm
{binary}}^{*}(\alpha)}$, then all tests are asymptotically powerless.

\item[2.] Let\vspace*{1pt} $A=\sqrt{\frac{2t\log(p)}{r_*}}$. If $t>{\rho_{\mathrm
{binary}}^{*}(\alpha)}$, then the Higher Criticism Test is
asymptotically powerful.
\end{longlist}
\end{theorem}

%

\begin{rem}
The assumptions on the design matrix in Theorem \ref
{theoremwalargersparse} is weaker than the assumptions in Theorem
\ref{theoremwasmallr}. In particular, one is allowed to go beyond
$\llvert \bigcup_{i \notin\Omega^*} S_i\rrvert \ll p$ in Theorem \ref
{theoremwasmallr} as long as the condition \textup{(C3)} is satisfied
with $\gamma=\log(p)$. This is expected since the conditions under
which all tests are asymptotically powerless irrespective of sample
size are often more stringent.
\end{rem}

%
\begin{rem}
Theorem \ref{theoremwalargersparse} states that the Higher
Criticism Test attains the sharp detection boundary in the sparse
regime. Note that the difference in the denominators of $A$ in the
statement of upper and lower bound in Theorem \ref
{theoremwalargersparse} is unavoidable and the difference vanishes
asymptotically if $r^*/r_*\rightarrow1$. This is expected since the
detection boundary depends on the column norms of the design matrix.
\end{rem}


\section{Simulation studies}\label{sec8}\label{secNE}
We complement our study with some numerical simulations which
illustrate the empirical performance of the test statistics described
in earlier sections for finite sample sizes. Since detection complexity
of the general weakly correlated binary design matrices depend on the
behavior of $\operatorname{ANOVA}(r)$ type designs, we only provide
simulations for strong
one-way ANOVA type design. Let $X$ be a balanced design matrix with
$p={}$10,000 covariates and $r$ replicates per covariate. For different
values of sparsity index $\alpha\in(0,1)$ and $r$, we study the
performance of Higher Criticism Test, GLRT and Max Test, respectively,
for different values of $t$, where $t$ which corresponds to $A=\sqrt
{\frac{2(\rho^{*}_{\mathrm{logistic}}(\alpha)+t)\log(p)}{r}}$.

Following \citet{Candes}, the performance of each of the three methods
is computed in terms of the empirical risk defined as the sum of
probabilities of type I and II errors achievable across all thresholds.
The errors are averaged over $300$ trials.
Even though the theoretical calculation of null distribution of the
Higher Criticism Test statistic computed from $p$-values remains a
challenge, we found that using
the $p$-value based statistic $\max_{1\leq j \leq p/2}\sqrt
{p}\frac{(j/p)-q_{(j)}}{\sqrt{q_{(j)}(1-q_{(j)})}}$
yielded expected results similar to our version of discretized Higher Criticism.

To be precise,
the performance of the test based on\break $\max_{1\leq j \leq
{p}/{2}}\sqrt{p}\frac{({j}/{p})-q_{(j)}}{\sqrt
{q_{(j)}(1-q_{(j)})}}$ was similar to the performance of the test based\vspace*{1pt}
on $\thc$ defined in Section \ref{sec5.2}. Note that this statistic is
different from $\thc'$ in that the maximum is taken over the
first $\frac{p}{2}$ elements instead of all $p$ of them. The main
reason for this is the fact that, as noted by \citet{Jin1}, the
information about the signal in the sample lies away from the extreme
$p$-values. The GLRT is based on $\tglrt$ as defined in Section \ref{sec5.1}
and the Max Test is based on the test statistic defined in Section~\ref{sec6.2.2.3}.

\begin{figure}

\includegraphics{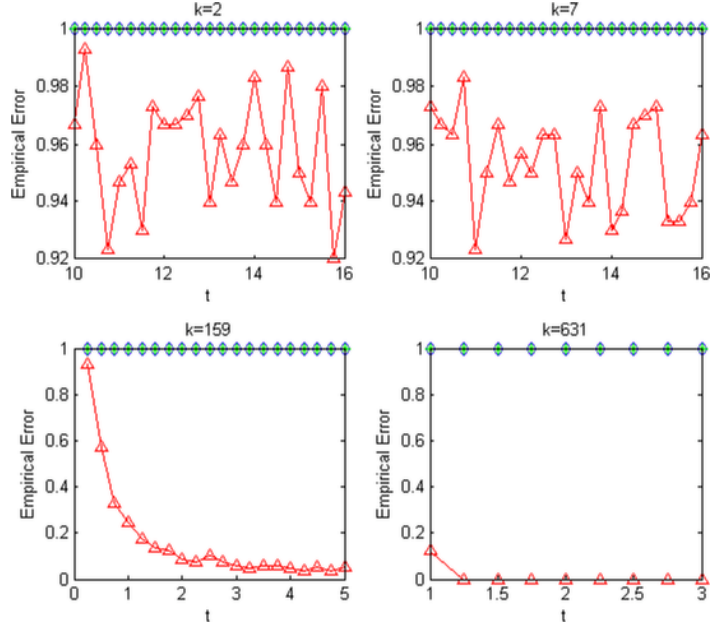}

\caption{Simulation results are for $p={}$10,000 and $r=\lceil\sqrt
{\log(p)}\rceil=4$. Sparsity level $k$ is indicated below each plot.
In each plot, the empirical risk of
each method [GLRT (triangles); Higher Criticism (diamonds);
Max Test (stars)] is
plotted against $t$ which corresponds to~$A=\sqrt{\frac{\max\{2(\rho
^{*}_{\mathrm{logistic}}(\alpha)+t),0\}\log(p)}{r}}$.} \label{fig3paper2}\vspace*{-3pt}
\end{figure}

\begin{figure}

\includegraphics{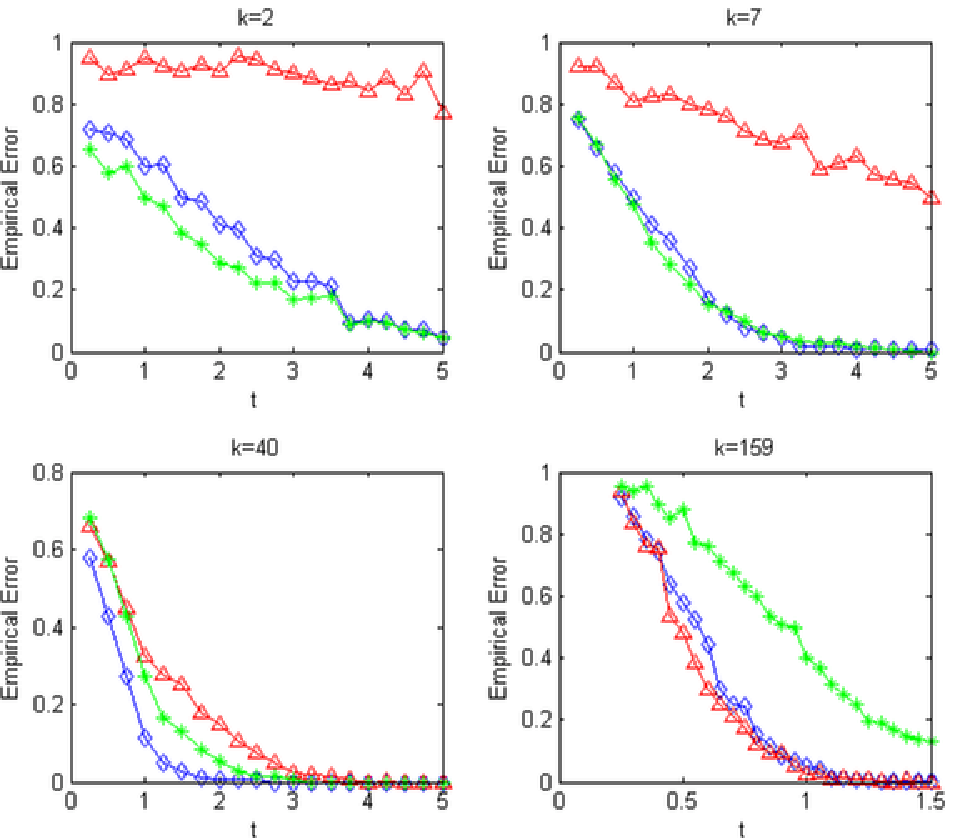}

\caption{Simulation results are for $p={}$10,000 and $r=\lceil(\log
(p))^5\rceil={}$66,280. Sparsity level $k$ is indicated below each plot.
In each plot, the empirical risk of
each method [GLRT (triangles); Higher Criticism (diamonds);
Max Test (stars)] is
plotted against $t$ which corresponds to~$A=\sqrt{\frac{2(\rho
^{*}_{\mathrm{logistic}}(\alpha)+t)\log(p)}{r}}$.} \label{fig2paper2}
\end{figure}

The results are reported in Figures~\ref{fig3paper2}~and~\ref
{fig2paper2}. For $r=\sqrt{\log(p)} \ll\log(p)$ and $k=2,7$ which
corresponds to $k \ll\sqrt{p}$, {that is}, the sparse regime,
we can see that all tests are asymptotically powerless in Figure~\ref
{fig3paper2} which is expected from the theoretical results. However,
even when $r=\lceil\sqrt{\log(p)}\rceil\ll\log(p)$, for the dense
regime, and $k=159$ and $631$, we see from Figure~\ref{fig3paper2}
that the GLRT is asymptotically powerful whereas the other two tests
are asymptotically powerless.
Once $r$ is much larger than $\log(p)$ in Figure~\ref{fig2paper2},
our observations are similar to \citet{Candes}. Here, we employ
simulations for $k=2,7,40$ which correspond to the sparse regime and
for $k=159$ which corresponds to the dense regime. We note that the
performance of GLRT improves very quickly as the sparsity decreases and
begins dominating the Max Test. 
The performance of the Max Test follows the opposite pattern with
errors of testing increasing as $k$ increases. The Higher Criticism
Test, however, continues to have good performance across the different
sparsity levels once $r \gg\log(p)$.

%

%
%

\section{Discussions}\label{sec9}
In this paper, we study testing of the global null hypothesis against
sparse alternatives in the context of general binary regression. We
show that, unlike Gaussian regression, the problem depends not only on
signal sparsity and strength, but also heavily on a design matrix
sparsity index. We
provide conditions on the design matrix which render all tests
asymptotically powerless irrespective of signal strength. In the
special case of design matrices with binary entries and certain
sparsity structures, we derive the lower and upper bounds for the
testing problem in both dense (rate optimal) and sparse regimes (sharp
including constants). In this context, we also develop a version of the
Higher Criticism Test statistic applicable for binary data which
attains the sharp detection boundary in the sparse regime.

In this paper, we constructed tests by combining tests based on
$Z$-statistics from the orthogonal part and the nonorthogonal part of the
$\bX$. In particular, we combine procedures based on $Z_j$ and
$Z_j^{\mathbf{G}}$ separately. This helps us achieve optimal rates for
upper bounds on testing errors under the same conditions required for
lower bounds in these problems. Indeed, one can consider constructing
GLRT and Higher Criticism Test using $Z$-statistics constructed based on
whole $\bX$, {that is}, based on $Z_j^{\bX}=X_j^T y, j=1,\ldots,p$
directly. We could obtain similar results based on the
combined $Z$-statistics under stronger structural information on $\mathbf
{G}$ than what we require here.

In particular, the conditions regarding the relative size of $\mathbf
{G}$ with respect to the orthogonal part of the design matrix, can be
substantially relaxed if more structural assumptions on $\mathbf{G}$
are made. For example,
for sequencing data, as observed in the Dallas Heart study data, for
people having more than one mutation, the locations of the mutations
are in fact usually clustered, due to linkage disequilibrium. For such
structures, strong results can be obtained. We omit those results here
due to space limitation. Future research is also needed to study the
detection boundary\vadjust{\goodbreak} for binary regression for more general design
matrices.

The study of detection boundaries associated with binary regression
models for a general design matrix is much more delicate. We allow in
this paper for a more general sparse design when the nonorthogonal
columns of the design matrix are sufficiently sparse and the number of
subjects with multiple nonzero entries in the design matrix are not too
large. Future research is needed to extend the results to a general
design matrix allowing a stronger correlation among the covariates~$X_j$'s.

\section*{Acknowledgements}
We would like to thank the Editor, Dr. Runze Li, the Associate Editors
and the reviewers for several
insightful comments which helped us improve the paper. Natesh S. Pillai
gratefully acknowledges the support from ONR.

\begin{supplement}[id=suppA]
\stitle{Supplement to ``Hypothesis testing for high-dimensional sparse binary regression''}
\slink[doi]{10.1214/14-AOS1279SUPP} 
\sdatatype{.pdf}
\sfilename{aos1279\_supp.pdf}
\sdescription{The supplementary material contain the proofs of all
theorems, propositions and supporting lemmas.}
\end{supplement}

%

\printaddresses

\begin{thebibliography}{18}
\bibitem[\protect\citeauthoryear{Arias-Castro, Cand{\`e}s and
Plan}{2011}]{Candes}
%
\begin{barticle}[mr]
\bauthor{\bsnm{Arias-Castro},~\bfnm{Ery}\binits{E.}},
\bauthor{\bsnm{Cand{\`e}s},~\bfnm{Emmanuel~J.}\binits{E.~J.}} \AND
\bauthor{\bsnm{Plan},~\bfnm{Yaniv}\binits{Y.}}
(\byear{2011}).
\btitle{Global testing under sparse alternatives: ANOVA, multiple
comparisons and the higher criticism}.
\bjournal{Ann. Statist.}
\bvolume{39}
\bpages{2533--2556}.
\bid{doi={10.1214/11-AOS910}, issn={0090-5364}, mr={2906877}}
\end{barticle}
%
\bptok{imsref}%
\endbibitem

\bibitem[\protect\citeauthoryear{Baraud}{2002}]{baraud1}
%
\begin{barticle}[mr]
\bauthor{\bsnm{Baraud},~\bfnm{Yannick}\binits{Y.}}
(\byear{2002}).
\btitle{Non-asymptotic minimax rates of testing in signal detection}.
\bjournal{Bernoulli}
\bvolume{8}
\bpages{577--606}.
\bid{issn={1350-7265}, mr={1935648}}
\end{barticle}
%
\bptok{imsref}%
\endbibitem

\bibitem[\protect\citeauthoryear{Cai, Jeng and Jin}{2011}]{Cai1}
%
\begin{barticle}[mr]
\bauthor{\bsnm{Cai},~\bfnm{T.~Tony}\binits{T.~T.}},
\bauthor{\bsnm{Jeng},~\bfnm{X.~Jessie}\binits{X.~J.}} \AND
\bauthor{\bsnm{Jin},~\bfnm{Jiashun}\binits{J.}}
(\byear{2011}).
\btitle{Optimal detection of heterogeneous and heteroscedastic mixtures}.
\bjournal{J. R. Stat. Soc. Ser. B Stat. Methodol.}
\bvolume{73}
\bpages{629--662}.
\bid{doi={10.1111/j.1467-9868.2011.00778.x}, issn={1369-7412}, mr={2867452}}
\end{barticle}
%
\bptok{imsref}%
\endbibitem

\bibitem[\protect\citeauthoryear{1000 Genomes Project
Consortium}{2012}]{10002012integrated}
%
\begin{barticle}[author]
\bauthor{\bsnm{1000 Genomes Project Consortium} and others}
(\byear{2012}).
\btitle{An integrated map of genetic variation from 1,092 human genomes}.
\bjournal{Nature}
\bvolume{491}
\bpages{56--65}.
\end{barticle}
%
\bptok{imsref}%
\endbibitem

\bibitem[\protect\citeauthoryear{Donoho and Jin}{2004}]{Jin1}
%
\begin{barticle}[mr]
\bauthor{\bsnm{Donoho},~\bfnm{David}\binits{D.}} \AND
\bauthor{\bsnm{Jin},~\bfnm{Jiashun}\binits{J.}}
(\byear{2004}).
\btitle{Higher criticism for detecting sparse heterogeneous mixtures}.
\bjournal{Ann. Statist.}
\bvolume{32}
\bpages{962--994}.
\bid{doi={10.1214/009053604000000265}, issn={0090-5364}, mr={2065195}}
\end{barticle}
%
\bptok{imsref}%
\endbibitem

\bibitem[\protect\citeauthoryear{Fu et~al.}{2013}]{fu2013analysis}
%
\begin{barticle}[author]
\bauthor{\bsnm{Fu},~\bfnm{Wenqing}\binits{W.}},
\bauthor{\bsnm{O'Connor},~\bfnm{Timothy~D.}\binits{T.~D.}},
\bauthor{\bsnm{Jun},~\bfnm{Goo}\binits{G.}},
\bauthor{\bsnm{Kang},~\bfnm{Hyun~Min}\binits{H.~M.}},
\bauthor{\bsnm{Abecasis},~\bfnm{Goncalo}\binits{G.}},
\bauthor{\bsnm{Leal},~\bfnm{Suzanne~M.}\binits{S.~M.}},
\bauthor{\bsnm{Gabriel},~\bfnm{Stacey}\binits{S.}},
\bauthor{\bsnm{Rieder},~\bfnm{Mark~J.}\binits{M.~J.}},
\bauthor{\bsnm{Altshuler},~\bfnm{David}\binits{D.}},
\bauthor{\bsnm{Shendure},~\bfnm{Jay}\binits{J.}} \betal{et~al.}
(\byear{2013}).
\btitle{Analysis of 6,515 exomes reveals the recent origin of most
human protein-coding variants}.
\bjournal{Nature}
\bvolume{493}
\bpages{216--220}.
\end{barticle}
%
\bptok{imsref}%
\endbibitem

\bibitem[\protect\citeauthoryear{Hall and Jin}{2010}]{Jin2}
%
\begin{barticle}[mr]
\bauthor{\bsnm{Hall},~\bfnm{Peter}\binits{P.}} \AND
\bauthor{\bsnm{Jin},~\bfnm{Jiashun}\binits{J.}}
(\byear{2010}).
\btitle{Innovated higher criticism for detecting sparse signals in
correlated noise}.
\bjournal{Ann. Statist.}
\bvolume{38}
\bpages{1686--1732}.
\bid{doi={10.1214/09-AOS764}, issn={0090-5364}, mr={2662357}}
\end{barticle}
\bptok{imsref}%
\endbibitem

\bibitem[\protect\citeauthoryear{Ingster and Suslina}{2003}]{Ingster4}
%
\begin{bbook}[mr]
\bauthor{\bsnm{Ingster},~\bfnm{Yu.~I.}\binits{Yu.~I.}} \AND
\bauthor{\bsnm{Suslina},~\bfnm{I.~A.}\binits{I.~A.}}
(\byear{2003}).
\btitle{Nonparametric Goodness-of-Fit Testing Under {G}aussian Models}.
\bseries{Lecture Notes in Statistics}
\bvolume{169}.
\bpublisher{Springer},
\blocation{New York}.
\bid{doi={10.1007/978-0-387-21580-8}, mr={1991446}}
\end{bbook}
%
\bptok{imsref}%
\endbibitem

\bibitem[\protect\citeauthoryear{Ingster, Tsybakov and
Verzelen}{2010}]{Ingster5}
%
\begin{barticle}[mr]
\bauthor{\bsnm{Ingster},~\bfnm{Yuri~I.}\binits{Y.~I.}},
\bauthor{\bsnm{Tsybakov},~\bfnm{Alexandre~B.}\binits{A.~B.}} \AND
\bauthor{\bsnm{Verzelen},~\bfnm{Nicolas}\binits{N.}}
(\byear{2010}).
\btitle{Detection boundary in sparse regression}.
\bjournal{Electron. J. Stat.}
\bvolume{4}
\bpages{1476--1526}.
\bid{doi={10.1214/10-EJS589}, issn={1935-7524}, mr={2747131}}
\end{barticle}
%
\bptok{imsref}%
\endbibitem

\bibitem[\protect\citeauthoryear{Koml{\'o}s, Major and Tusn{\'a}dy}{1975}]{kmt}
%
\begin{barticle}[mr]
\bauthor{\bsnm{Koml{\'o}s},~\bfnm{J.}\binits{J.}},
\bauthor{\bsnm{Major},~\bfnm{P.}\binits{P.}} \AND
\bauthor{\bsnm{Tusn{\'a}dy},~\bfnm{G.}\binits{G.}}
(\byear{1975}).
\btitle{An approximation of partial sums of independent {${\rm RV}$}'s
and the sample {${\rm DF}$}. {I}}.
\bjournal{Z. Wahrsch. Verw. Gebiete}
\bvolume{32}
\bpages{111--131}.
\bid{mr={0375412}}
\end{barticle}
%
\bptok{imsref}%
\endbibitem

\bibitem[\protect\citeauthoryear{Lee et~al.}{2014}]{lee2014}
%
\begin{barticle}[author]
\bauthor{\bsnm{Lee},~\bfnm{S.}\binits{S.}},
\bauthor{\bsnm{Abecasis},~\bfnm{G.}\binits{G.}},
\bauthor{\bsnm{Boehnke},~\bfnm{M.}\binits{M.}} \AND
\bauthor{\bsnm{Lin},~\bfnm{X.}\binits{X.}}
(\byear{2014}).
\btitle{Analysis of rare variants in sequencing-based association studies}.
\bjournal{The American Journal of Human Genetics}
\bvolume{95}
\bpages{5--23}.
\end{barticle}
%
\bptok{imsref}%
\endbibitem

\bibitem[\protect\citeauthoryear{Mukherjee, Pillai and
Lin}{2014}]{MukherjeePillaiLin2014}
%
\begin{bmisc}[author]
{\bauthor{\bsnm{Mukherjee},~\binits{R.}},
\bauthor{\bsnm{Pillai},~\binits{N. S.}} \AND
\bauthor{\bsnm{Lin},~\binits{X.}}}
(\byear{2014}).
\bhowpublished{Supplement to ``Hypothesis testing for high-dimensional
sparse binary regression.''
DOI:\doiurl{10.1214/14-AOS1279SUPP}}.
\bptok{imsref}%
\end{bmisc}
%
\bptok{imsref}%
\endbibitem

\bibitem[\protect\citeauthoryear{Nelson et~al.}{2012}]{nelson2012abundance}
%
\begin{barticle}[author]
\bauthor{\bsnm{Nelson},~\bfnm{Matthew~R.}\binits{M.~R.}},
\bauthor{\bsnm{Wegmann},~\bfnm{Daniel}\binits{D.}},
\bauthor{\bsnm{Ehm},~\bfnm{Margaret~G.}\binits{M.~G.}},
\bauthor{\bsnm{Kessner},~\bfnm{Darren}\binits{D.}},
\bauthor{\bsnm{Jean},~\bfnm{Pamela~St}\binits{P.~S.}},
\bauthor{\bsnm{Verzilli},~\bfnm{Claudio}\binits{C.}},
\bauthor{\bsnm{Shen},~\bfnm{Judong}\binits{J.}},
\bauthor{\bsnm{Tang},~\bfnm{Zhengzheng}\binits{Z.}},
\bauthor{\bsnm{Bacanu},~\bfnm{Silviu-Alin}\binits{S.-A.}},
\bauthor{\bsnm{Fraser},~\bfnm{Dana}\binits{D.}} \betal{et~al.}
(\byear{2012}).
\btitle{An abundance of rare functional variants in 202 drug target
genes sequenced in 14,002 people}.
\bjournal{Science}
\bvolume{337}
\bpages{100--104}.
\end{barticle}
\bptok{imsref}%
\endbibitem

\bibitem[\protect\citeauthoryear{Plan and Vershynin}{2013a}]{plan1}
%
\begin{barticle}[mr]
\bauthor{\bsnm{Plan},~\bfnm{Yaniv}\binits{Y.}} \AND
\bauthor{\bsnm{Vershynin},~\bfnm{Roman}\binits{R.}}
(\byear{2013}a).
\btitle{One-bit compressed sensing by linear programming}.
\bjournal{Comm. Pure Appl. Math.}
\bvolume{66}
\bpages{1275--1297}.
\bid{doi={10.1002/cpa.21442}, issn={0010-3640}, mr={3069959}}
\end{barticle}
%
\bptok{imsref}%
\endbibitem

\bibitem[\protect\citeauthoryear{Plan and Vershynin}{2013b}]{plan2}
%
\begin{barticle}[mr]
\bauthor{\bsnm{Plan},~\bfnm{Yaniv}\binits{Y.}} \AND
\bauthor{\bsnm{Vershynin},~\bfnm{Roman}\binits{R.}}
(\byear{2013}b).
\btitle{Robust 1-bit compressed sensing and sparse logistic regression:
A convex programming approach}.
\bjournal{IEEE Trans. Inform. Theory}
\bvolume{59}
\bpages{482--494}.
\bid{doi={10.1109/TIT.2012.2207945}, issn={0018-9448}, mr={3008160}}
\end{barticle}
%
\bptok{imsref}%
\endbibitem

\bibitem[\protect\citeauthoryear{Tang et~al.}{2014}]{tang2013large}
%
\begin{barticle}[author]
\bauthor{\bsnm{Tang},~\bfnm{Huayang}\binits{H.}},
\bauthor{\bsnm{Jin},~\bfnm{Xin}\binits{X.}},
\bauthor{\bsnm{Li},~\bfnm{Yang}\binits{Y.}},
\bauthor{\bsnm{Jiang},~\bfnm{Hui}\binits{H.}},
\bauthor{\bsnm{Tang},~\bfnm{Xianfa}\binits{X.}},
\bauthor{\bsnm{Yang},~\bfnm{Xu}\binits{X.}},
\bauthor{\bsnm{Cheng},~\bfnm{Hui}\binits{H.}},
\bauthor{\bsnm{Qiu},~\bfnm{Ying}\binits{Y.}},
\bauthor{\bsnm{Chen},~\bfnm{Gang}\binits{G.}},
\bauthor{\bsnm{Mei},~\bfnm{Junpu}\binits{J.}} \betal{et~al.}
(\byear{2014}).
\btitle{A large-scale screen for coding variants predisposing to psoriasis}.
\bjournal{Nature Genetics}
\bvolume{46}
\bpages{40--50}.
\end{barticle}
\bptok{imsref}%
\endbibitem

\bibitem[\protect\citeauthoryear{Victor et~al.}{2004}]{dallas}
%
\begin{barticle}[author]
\bauthor{\bsnm{Victor},~\bfnm{Ronald~G.}\binits{R.~G.}},
\bauthor{\bsnm{Haley},~\bfnm{Robert~W.}\binits{R.~W.}},
\bauthor{\bsnm{Willett},~\bfnm{DuWayne~L.}\binits{D.~L.}},
\bauthor{\bsnm{Peshock},~\bfnm{Ronald~M.}\binits{R.~M.}},
\bauthor{\bsnm{Vaeth},~\bfnm{Patrice~C.}\binits{P.~C.}},
\bauthor{\bsnm{Leonard},~\bfnm{David}\binits{D.}},
\bauthor{\bsnm{Basit},~\bfnm{Mujeeb}\binits{M.}},
\bauthor{\bsnm{Cooper},~\bfnm{Richard~S.}\binits{R.~S.}},
\bauthor{\bsnm{Iannacchione},~\bfnm{Vincent~G.}\binits{V.~G.}},
\bauthor{\bsnm{Visscher},~\bfnm{Wendy~A.}\binits{W.~A.}} \betal{et~al.}
(\byear{2004}).
\btitle{The Dallas Heart Study: A population-based probability sample
for the multidisciplinary study of ethnic differences in cardiovascular health}.
\bjournal{The American Journal of Cardiology}
\bvolume{93}
\bpages{1473--1480}.
\end{barticle}
\bptok{imsref}%
\endbibitem

\bibitem[\protect\citeauthoryear{Wald}{1950}]{wald1}
%
\begin{bbook}[author]
\bauthor{\bsnm{Wald},~\bfnm{Abraham}\binits{A.}}
(\byear{1950}).
\btitle{Statistical Decision Functions.}
\bpublisher{Chelsea},
\blocation{New York}.
\end{bbook}
\bptok{imsref}%
\endbibitem
\end{thebibliography}
\end{document}